\documentclass[10pt]{article}
\usepackage{psfrag,
graphicx, harvard,
epsf,
euscript,amsfonts,amsmath,latexsym,amssymb,mathrsfs
}
\usepackage[active]{srcltx}
\newtheorem{lemma}{Lemma}
\newtheorem{theorem}{Theorem}
\newtheorem{proposition}{Proposition}

\newtheorem{remark}{Remark}
\newtheorem{definition}{Definition}

\newcommand{\cK}{{\cal K}}

\newcommand{\cR}{{\cal R}}

\newcommand{\bL}{{\mathbb L}}

\newcommand{\bN}{{\mathbb N}}

\newcommand{\bR}{{\mathbb R}}

\newcommand{\bZ}{{\mathbb Z}}
\newcommand{\sC}{{\mathscr C}}

\newcommand{\sG}{{\mathscr G}}
\newcommand{\sH}{{\mathscr H}}

\newcommand{\sM}{{\mathscr M}}

\newcommand{\sP}{{\mathscr P}}

\newcommand{\rP}{\mathrm{P}}
\newcommand{\rE}{\mathrm{E}}

\newcommand{\rZ}{\mathrm{Z}}

\newcommand{\rd}{\mathrm{d}}

\renewcommand{\kappa}{\varkappa}
\newcommand{\pr} {\par \noindent{\bf Proof\,:~}}
\newcommand{\epr}{\hfill\hbox{\hskip 4pt
                \vrule width 5pt height 6pt depth 1.5pt}\vspace{0.5cm}\par}
\usepackage[colorlinks=true, linkcolor=blue, pdfborder={0 0 0}]{hyperref}
\begin{document}
\author{
A.~Goldenshluger\thanks{The author is grateful to Gideon Weiss for attracting his
attention  to the problem studied in this paper, and to Oleg Lepski 
for 
useful discussions and suggestions. Part of this work has been done  while
the author was visiting NYU Shanghai.}
\\*[2mm]
{\small
Department of Statistics}
\\
{\small
University of Haifa
}
\\
{\small Haifa 31905, Israel}
}
\title{Nonparametric estimation of service time distribution 
in 
 the $M/G/\infty$ queue and related estimation
problems}
\date{}
\maketitle
\begin{abstract} 
The subject of this paper is the problem of estimating service time distribution
of the $M/G/\infty$ queue from incomplete data on the queue. 
The goal is to estimate $G$ from 
observations of the queue--length process at the points of  the regular grid
on a fixed time interval. We propose an estimator and analyze its accuracy 
over a family of target service time distributions. 
The original $M/G/\infty$ problem is closely related to 
the problem of estimating derivatives of
the covariance function of a  stationary Gaussian process.
We consider the latter problem and
 derive lower bounds on the minimax risk.
The obtained results strongly suggest that the proposed estimator of the service time
distribution is rate optimal.
\end{abstract}

\vspace*{1em} \noindent {\bf Short Title:} Nonparametric estimation in the $M/G/\infty$ queue


\noindent {\bf 2000 AMS Subject Classification} : 
\section{Introduction}
Suppose that customers arrive 
at a system at time instances $\{\tau_j, j\in \bZ\}$, 
obtain service upon arrival, and leave the system at time instances $\{y_j, j\in \bZ\}$
after the service is completed. 
A $j$th customer arriving at  $\tau_j$ requires service time $\sigma_j$, so
that its departure epoch 
is $y_j=\tau_j+\sigma_j$.
If $\{\tau_j, j\in\bZ\}$  is a realization of a stationary Poisson process
on $\bR$, and   
$\{\sigma_j, j\in\bZ\}$ are non--negative independent  
random variables with common distribution $G$, independent of  
$\{\tau_j, j\in\bZ\}$, then the above description corresponds 
to the $M/G/\infty$ queueing system. 
In this paper we are interested in estimating
service time distribution $G$ from incomplete data on the queue.
\par 
The $M/G/\infty$ system is perhaps one of the most widely studied models in queueing theory;  
its probabilistic properties are fairly well understood. 
However statistical inference in 
such models has attracted  little attention.
\par 
The problem of estimating service time distribution $G$ in 
the $M/G/\infty$ queue has been studied
under different assumptions on the 
available data. 
The following three observation schemes
have been considered in the literature: 
\begin{itemize}
\item[(i)] observation of arrival $\{\tau_j, j\in \bZ\}$ and 
departure $\{y_j, j\in \bZ\}$ epochs without their matchings;  
\item[(ii)] observation of the queue--length
(number--of--busy--servers) process $\{X(t)\}$; 
\item[(iii)] observation of the busy--period process 
\mbox{$\{{\bf 1}(X(t)>0)\}$}.
\end{itemize}
\par 
We note that observation schemes (i) and (ii) are equivalent up to 
initial conditions on the queue length. 
In particular, arrival and departure epochs are uniquely determined 
by the queue--length process, 
while the queue length can be reconstructed 
from the input--output data provided that the initial state of the queue is known.
\par 
In setting (i)
\citeasnoun{Brown} 
proposed an estimator of $G$ which is based on the idea 
of pairing every departure epoch with the closest arrival epoch to the left.
Differences between 
these epochs constitute an ergodic stationary random 
sequence whose marginal distribution is related to the service time distribution $G$ by a simple formula.   
Then estimation of $G$ can be achieved by inverting the formula
and substituting the  empirical marginal distribution of the differences. 
\citeasnoun{Brown} proved that the proposed  estimator is   consistent.
 Recently \citeasnoun{Weiss} extended the work of Brown; 
 they showed  
that pairing of a departure 
epoch with the $r$--closest arrival epoch to the left can be worthwhile. 
\par 
Nonparametric estimation of service time distribution $G$ under observation schemes (ii) and (iii)
was considered in \citeasnoun{Bingham}. It is well known that in the steady state the queue--length process
$\{X(t)\}$ is stationary with Poisson marginal distribution and correlation function 
\begin{equation}\label{eq:correlation}
H(t)=1-G^*(t),\;\;\;
G^*(t):=\Big[\int_0^\infty [1-G(x)]\rd x\Big]^{-1} \int_0^t [1-G(x)]\rd x;
\end{equation}
see, e.g., \citeasnoun{Benes} and \citeasnoun{Reynolds}.
This fact suggests that function $G^*$ can be reconstructed by estimating 
correlation function of the queue--length process. The work of \citeasnoun{Bingham} 
discusses this approach and 
provides standard
results from the time series literature for estimators of $G^*$. 
The idea of reconstructing  the service time 
distribution from correlation structure of the queue--length process was also exploited by
\citeasnoun{Pickands}. The model  considered in that paper assumes that 
a Poisson number
of customers arrives at discrete times $1,2,\ldots, T$, and service times are i.i.d.
random variables taking values in the set of non--negative integer numbers. 
In this discrete setting  estimation of the service time distribution is equivalent 
to estimating a linear form of the correlation function of the queue--length process. 
For the latter problem standard results from the time series literature are applicable.
Other related work is reported in 
\citeasnoun{Brillinger}, \citeasnoun{Bingham-Dunham}, \citeasnoun{Hall},
\citeasnoun{Moulines},
\citeasnoun{Grubel}, \citeasnoun{cornelia}; see \citeasnoun{Weiss} for additional references.  
\par
Although  estimation of $G$ under different observation schemes
was considered in the literature, the most interesting and important statistical questions
remain to be open.  In particular, it is not clear 
what is the achievable
estimation accuracy in such problems, and  how to construct optimal estimators.
The goal of  this paper is to shed light on some of these issues. 
\par
In this work we adopt minimax approach for measuring estimation accuracy. It is assumed that 
the estimated distribution $G$ belongs
to a given functional class, and accuracy of any estimator is measured by its worst--case 
mean squared error on the class.  The functional class is defined in terms of   
restrictions on  smoothness and tail behavior of $G$ (for precise definitions see Section~\ref{sec:problem}). 
We concentrate on the observation scheme~(ii) when the queue--length process is observed on 
a fixed interval at the points of
the regular grid. 
We want to estimate $G$ at a fixed point using such observations. 
From now on we will refer to this setting as 
{\em the $M/G/\infty$ estimation problem}.
\par 
We develop an estimator of $G$ which is based 
on the relationship between distribution $G$ and 
covariance function of the 
queue--length process, as discussed in \citeasnoun{Bingham}
and \citeasnoun{Pickands} [cf. (\ref{eq:correlation})].  In particular, 
estimating $G$ at a fixed point is reduced to estimating
derivative of the covariance function of the queue--length process at this point.  
We analyze accuracy of our estimator over a suitable class of target distributions and derive 
an upper bound
on the maximal risk. The upper bound is expressed in terms of  the functional class parameters and 
the observation horizon. The problem of estimating the arrival rate is discussed as well.
\par 
A natural question is:  what is the achievable estimation accuracy in the $M/G/\infty$
problem? This question calls for a lower bound on the minimax risk.
Since explicit formulas for finite dimensional distributions of the queue--length process
in the $M/G/\infty$ model 
are not available, derivation of lower bounds on the minimax risk 
seems to be analytically intractable. Therefore, driven by a 
Gaussian approximation to the queue--length process,
we consider a closely related estimation problem 
for a Gaussian model. Specifically, 
let $\{X(t), t\in \bR\}$ be a continuous--time stationary Gaussian
process  which is observed at the points of a regular grid on a given time interval.
Using such discrete observations we want to estimate the derivative of the covariance function
of $\{X(t), t\in \bR\}$. We derive a lower bound on the minimax risk in this problem, and show
that under suitable conditions it converges to zero at the same rate as 
the risk of our estimator in the $M/G/\infty$
estimation problem. This fact strongly suggests that our estimator
of the service time distribution
is rate--optimal.  
\par 
The problem of  estimating derivatives  
of covariance functions at a fixed point 
(or, more generally, linear functionals of covariance functions/spectral densities)
from discrete observations 
is interesting in its own right.  
Although various settings 
were considered in the literature, we are not aware of any work dealing with 
estimation of covariance function derivatives.
For discrete--time stationary processes asymptotic efficient
estimators of smooth functionals of the spectral density 
were proposed in \citeasnoun{Hasminskii};
see also 
\citeasnoun{Ginovyan}, where  continuous--time stationary
processes and continuous observations were considered. 
Nonparametric estimation of 
covariance functions for continuous--time stationary processes 
from
discrete observations is discussed in 
\citeasnoun{Hall94a} and \citeasnoun{Hall94b}.
For other related work we
refer to   \citeasnoun{Masry}, \citeasnoun{Haberzettl}, 
\citeasnoun{Srivastava} and references therein. 
\par 
 The rest of this paper is structured as follows. 
Section~\ref{sec:problem} contains formal statement of the $M/G/\infty$ estimation problem.
Section~\ref{sec:preliminaries}
presents some results on properties of the queue--length process; these results are instrumental
for subsequent
developments in the paper. 
In Section~\ref{sec:estimator} we consider the $M/G/\infty$ estimation problem,  
define our estimator and establish upper bounds on its maximal risk. 
Section~\ref{sec:arrival-rate}
deals with the problem of estimating the arrival rate in the $M/G/\infty$ queue.  
In Section~\ref{sec:covariance} we relate the $M/G/\infty$ problem to
the problem of estimating derivative of covariance
function of  a continuous--time stationary Gaussian process, and 
derive a lower bound on the minimax risk for the latter problem.
Proofs are given in Section~\ref{sec:proofs}.
\section{Problem formulation}\label{sec:problem}
Let $\{\tau_j, j\in \bZ\}$ be arrival epochs constituting a realization of
stationary Poisson process point process of intensity $\lambda$ on the real line. The service times
$\{\sigma_j, j\in \bZ\}$ are 
positive  independent random variables with common distribution $G$, independent of
$\{\tau_j, j\in\bZ\}$. 
Assume that 
the system is in the steady state; then the queue--length
process $\{X(t), t\in \bR\}$ is given by 
\begin{equation}\label{eq:X}
 X(t) = \sum_{j\in \bZ} {\bf 1}(\tau_j \leq t, \sigma_j > t-\tau_j),\;\;\;t\in\bR.
\end{equation}
Suppose that $X(t)$ is observed on the time interval $[0, T]$ at the points of the regular
grid $t_i=i\delta$, $i=1,\ldots, n$, where $\delta>0$ is the sampling interval, and $T=n\delta$. Denote 
$X^n=(X(t_1),\ldots, X(t_n))\in \bR_+^n$. 
Our goal is to estimate the distribution function $G$ at single given point $x_0\in \bR_+$
using observation $X^n$. In Section~\ref{sec:arrival-rate}
we also discuss the problem of estimating the arrival rate $\lambda$
from observation $X^n$.
\par
Distribution of the observation $X^n$  is fully characterized by the service time distribution $G$ and
by the arrival rate $\lambda$.
From now on $\rP_{G,\lambda}$ stands for the probability measure 
generated by $\{\tau_j, j\in \bZ\}$ and $\{\sigma_j, j\in \bZ\}$
when $\sigma_j$'s are distributed $G$, and the arrival rate is $\lambda$. Correspondingly, $\rE_{G,\lambda}$ is
the expectation with respect to $\rP_{G, \lambda}$. In the problem of estimating $G$ when the arrival 
rate $\lambda$ is known, 
we use notation $\rP_G$ and $\rE_G$ for the probability measure and expectation respectively.
\par 
By  estimator $\hat{G}(x_0)=\hat{G}(X^n; x_0)$ of $G(x_0)$ we mean any measurable function
of the observation $X^n$. We adopt minimax approach for measuring estimation accuracy.
Let $\sG$ be a class of distribution functions; then accuracy of $\hat{G}(x_0)$ is measured
by the maximal mean squared risk over the class:
\[
 \cR_{x_0}[\hat{G}; \sG]= \sup_{G\in \sG} 
\Big[\rE_G\, \big|\hat{G}(x_0) - G(x_0)\big|^2\Big]^{1/2}.
\]
The minimax risk is defined by 
\[
\cR_{x_0}^*[\sG]=\inf_{\hat{G}} \cR_{x_0}[\hat{G}; \sG], 
\]
where $\inf$ is taken over all possible estimators.   We want to develop  a {\em rate--optimal}
({\em optimal in order}) estimator $\tilde{G}(x_0)$ such that
\[ 
\cR_{x_0}[\tilde{G};\sG] \leq C \cR_{x_0}^*[\sG], 
\]
where $C$ is a constant independent of the observation horizon  $T$ and 
the sampling interval $\delta$. 
\par 
In the problem of estimating the arrival rate $\lambda$ from observation $X^n$
the estimation accuracy is measured similarly. If $\hat{\lambda}=\hat{\lambda}(X^n)$ 
is an estimator 
of $\lambda$ (a measurable function of $X^n$) then the maximal risk of $\hat{\lambda}$ is defined by 
\[
 \cR[\hat{\lambda}; \sG] = \sup_{G\in \sG}\big[\rE_{G,\lambda} |\hat{\lambda} - \lambda|^2\big]^{1/2}.
\]
We will consider functional classes $\sG$ which impose restrictions on smoothness and tail behavior
of the distribution functions. The corresponding definitions are given in 
Section~\ref{sec:estimator}.
\section{Queue--length process}\label{sec:preliminaries}
Let  
\[
\tfrac{1}{\mu}:= 
\rE_G [\sigma] =\int_0^\infty [1-G(t)]\rd t < \infty
\]
with $\mu$ being {\em the service rate}, and 
let
$\rho:=\lambda/\mu$ be the {\em traffic intensity}. 
Define 
\begin{equation}\label{eq:H}
 H(t):=\mu \int_{t}^\infty [1-G(x)]\rd x= 
\Big[\int_0^\infty [1-G(x)]\rd x\Big]^{-1}\int_t^\infty [1-G(x)]\rd x, 
\;\;\;\;\; t\in \bR_+.
\end{equation}
\par  
The function $G^*:=1-H$ is often 
called {\em the stationary--excess cumulative distribution function} 
[see, e.g., \cite{Whitt85}]. If $G$ is a distribution function of an interval 
between points in a renewal process, then $G^*$  represents a distribution
function of the interval between arbitrary time  and the next renewal point.
In our context, the important role of $H$ stems from the fact that it is 
the correlation  function of the queue--length process $\{X(t), t\in \bR\}$; 
see Proposition~\ref{lem:queue-size}
below. 
\par 
Observe that $H(0)=1$, and $H$ is monotone decreasing on the positive real line. Although 
function $H$ is defined on $\bR_+$ only, it will be convenient to extend its definition
to the whole real line $\bR$ by setting $H(t)=H(-t)$ for $t<0$. 
From now on we use the suffix notation $X_i=X(t_i)=X(i\delta)$, $H_i=H(t_i)=H(i\delta)$, etc.
\begin{proposition}\label{lem:queue-size}
The following statements hold.
\begin{itemize}
 \item[{\rm (i)}] For any $t\in \bR$ the distribution of 
$X(t)$ is Poisson with parameter $\rho$.
\item[{\rm (ii)}] For any $t, s\in \bR$
\[
 \rE_{G} \big[X(t) X(s)\big] = \rho^2 + \rho H(t-s).
\]
\item[{\rm (iii)}] For any $\theta=(\theta_1,\ldots, \theta_n)$, $n\geq 1$, one has
\begin{eqnarray}
&&\log \rE_G \Big[\exp\Big\{\sum_{i=1}^n \theta_i X_i\Big\}\Big] \;=\; \rho S_n(\theta),
\label{eq:S-1}
\\ 
&& S_{n}(\theta) \;:=\; 
 \sum_{m=1}^{n} \big(e^{\theta_m}-1\big) +
 \sum_{k=1}^{n-1} H_k \sum_{m=k}^{n-1} 
 \big(e^{\theta_{m-k+1}}-1\big)e^{\sum_{i=m-k+2}^m\theta_i}
 \big(e^{\theta_{m+1}}-1\big).
\label{eq:S-2}
\end{eqnarray}
In particular, if $\theta^*:=(\vartheta,\ldots,\vartheta)$ for some $\vartheta\in \bR$ then   
\begin{eqnarray*} 
S_n(\theta^*) &=& n(e^\vartheta-1)
+ n(e^\vartheta -1)^2\sum_{k=1}^{n-1} \big(1-\tfrac{k}{n}\big)
e^{(k-1)\vartheta}H_k.
\end{eqnarray*}
\end{itemize}
\end{proposition}

\begin{remark}\mbox{}
\begin{enumerate}
 \item[{\rm (i)}] The statements (i) and (ii) are well known; in fact, they are  
immediate consequences of~(iii).
The first statement can be found in many textbooks  [see, e.g.,
\citeasnoun[p.~147]{parzen} and \citeasnoun[p.~19]{ross}], while the second one  
appears, e.g., in \citeasnoun{Benes} and \citeasnoun{Reynolds}.
As for the part~(iii),  \citeasnoun{lindley} considered the 
special case of $n=3$ and discussed heuristically
a derivation
for general $n$. However, we could not find 
formula (\ref{eq:S-1})--(\ref{eq:S-2})
in the literature, and, to the best of our knowledge, it is new.
This formula plays an important role in subsequent derivations. 
\item[{\rm (ii)}] The joint distribution of $X^n$ is the so--called multivariate Poisson; for details see, e.g., 
\citeasnoun[\S 2]{lindley} and \citeasnoun{Milne}. The statements~(i) and~(ii) show that $H$
is the correlation  function of the process $\{X(t), t\in\bR\}$.
\end{enumerate}
\end{remark}
\par 
It is instructive to realize  the form of
(\ref{eq:S-1})--(\ref{eq:S-2}) in the special case  $n=4$.  
Let $1\leq i\leq j\leq k\leq m\leq n$; then
\begin{eqnarray}
&& \tfrac{1}{\rho}
\log \rE_G \big[\exp\{\theta_1 X_i +\theta_2 X_j + \theta_3 X_k +\theta_4X_m\}\big]
\;=\;\sum_{l=1}^4 (e^{\theta_l}-1) 
\nonumber
\\*[2mm]
&&\;\;\;\;\;\; +\;H_{j-i}(e^{\theta_1}-1)(e^{\theta_2}-1) 
\;+\;
H_{k-i}(e^{\theta_1}-1)e^{\theta_2}(e^{\theta_3}-1)
\nonumber
\\*[2mm]
&&\;\;\;\;\;\;+\;H_{m-i}(e^{\theta_1}-1)e^{\theta_2+\theta_3}(e^{\theta_4}-1)
\;+\;H_{k-j}(e^{\theta_2}-1)(e^{\theta_3}-1) 
\label{eq:4-var}
\\*[2mm]
&&\;\;\;\;\;\;+\;
H_{m-j}(e^{\theta_2}-1) e^{\theta_3}(e^{\theta_4}-1)
\;+\;\;H_{m-k} (e^{\theta_3}-1)(e^{\theta_4}-1).
 \nonumber
\end{eqnarray}
As it is seen from the above formula, the first term on the 
right hand side of 
(\ref{eq:4-var}) coincides with the cumulant generating function of independent
Poisson random variables. The other terms are associated with   
all possible pairs of random variables. For every pair of random variables 
the corresponding term contains correlation between the variables, and 
factors $(1-e^\theta)$ and $e^\theta$, where 
$(1-e^{\theta})$--factors correspond to the pair, and  
$e^\theta$--factors correspond to the random variables ``sandwitched'' by the pair.    
\par 
The formula~(\ref{eq:4-var}) allows to compute mixed moments of the 
fourth order as presented in the next statement.
\begin{proposition}\label{cor:moments}
Let $1\leq i\leq  j\leq k\leq  m\leq n$; then 
\begin{eqnarray}
&& \rE_G \big[X_i X_j X_k X_m\big] \;=\; \rho^4 + 
\rho^3 \big( H_{j-i}+ H_{k-i} + H_{m-i}+H_{k-j}+H_{m-j} + H_{m-k}\big) 
\nonumber
\\*[2mm]
&&\;\;\;+\;\rho^2 \big(H_{k-i}+H_{m-j}+2H_{m-i}+H_{j-i}H_{m-k} + H_{k-i}H_{m-j} +H_{k-j}H_{m-i}\big)
+ \rho H_{m-i}.
\nonumber
\end{eqnarray}
More generally, for any $i,j,k,m \in \{1,\ldots, n\}$ and any subset $I$ of indexes 
$I\subseteq  \{i,j,k,m\}$
define 
$q_I=\max_{i_1,i_2\in I} |i_1-i_2|$.
Then  
\begin{eqnarray}
&& \rE_G [X_iX_jX_k] = \rho^3+\rho^2(H_{|i-j|}+H_{|k-j|}+H_{|k-i|})+ \rho H_{q_{\{i,j,k\}}}
\nonumber
\\*[2mm]
&& \rE_G \big[X_i X_j X_k X_m\big] 
\nonumber
\\*[2mm] 
&& \;\;=  \rho^4  + \rho^3 \big( H_{|j-i|}+H_{|k-i|}+H_{|m-i|}+H_{|k-j|}+H_{|m-j|} + H_{|m-k|}\big) + \rho H_{q_{\{i,j,k,m\}}}
\label{eq:moments-1}
\\*[2mm]
&&\;\; +\;\rho^2 \big[H_{q_{\{i,j,k\}}}+H_{q_{\{i,j,m\}}}+H_{q_{\{j,k,m\}}}+H_{q_{\{i,k,m\}}}+
H_{|j-i|}H_{|m-k|} + H_{|k-i|}H_{|m-j|} +H_{|k-j|}H_{|m-i|}\big].
\nonumber
\end{eqnarray}
\end{proposition}

\par
As a by--product of statement~(iii) in Proposition~\ref{lem:queue-size} we can easily obtain  
the following Gaussian approximation to 
finite dimensional distributions 
of the queue--length process $\{X(t), 0\leq t\leq T\}$.
\begin{proposition}\label{lem:Gaussian-approximation}
Consider a sequence of the 
$M/G/\infty$ queueing systems, 
$\{M_l/G/\infty, l=1,2,\ldots\}$, 
with the fixed service time distribution $G$ , and  with the $l$-th system characterized by 
the arrival rate $\lambda_l=l\lambda$, $\lambda>0$.
Let $X_l^n=(X_{l,1},\ldots, X_{l,n})=(X_l(t_1),\ldots,X_l(t_n))$ be the vector of observations 
of the queue--length
process (\ref{eq:X}) in the $l$-th system; then 
 \[
\frac{X_l^n - l \rho e_n}{\sqrt{l\rho}} \stackrel{d}{\to} {\cal N}_n \big(0, \Sigma(H)\big),\;\;\;l\to\infty,
\]
where $\rho=\lambda/\mu$, 
$e_n=(1,\ldots,1)\in \bR^n$, and $\Sigma(H):=\{H((i-j)\delta)\}_{i,j=1,\ldots, n}$.
\end{proposition}
\par 
The result of Proposition~\ref{lem:Gaussian-approximation} is well known; it is
in line with more general weak convergence results for queues
in \citeasnoun{Borovkov}, \citeasnoun{iglehart} and \citeasnoun{whitt}. 
The proof of Proposition~\ref{lem:Gaussian-approximation} follows immediately from 
Proposition~\ref{lem:queue-size}(iii), and it  is omitted.

\section{Estimation of service time distribution}\label{sec:estimator}
According to Proposition~\ref{lem:queue-size}(ii) the 
covariance function of the queue--length process is
\[
R(t):= {\rm cov}_G \{X(s), X(s+t)\}=\rho H(t).
\]
Therefore differentiation yields
\begin{equation}\label{eq:G-R-relation}
 1- G(t) = - \tfrac{1}{\lambda} R^\prime (t),\;\;\;t\in \bR_+.
\end{equation}
This relationship is the basis for construction of our estimator
of $G(x_0)$.
\subsection{Estimator construction}\label{sec:construction}
Let 
\[
\hat{\rho}_k= \tfrac{1}{n-k} \sum_{i=1}^{n-k} X_i,  \;\;\;k=0,1,\ldots, n-1,
\]
and  define
\begin{equation}\label{eq:R-estimator}
 \hat{R}_k = \tfrac{1}{n-k}\sum_{i=1}^{n-k} (X_i-\hat{\rho}_k) (X_{i+k}-\hat{\rho}_k), \;\;\;k=0,1,\ldots, n-1.
\end{equation}
Note that $\hat{R}_k$ is the empirical estimator 
of the covariance $R_k=R(k\delta)=\rho H(k\delta)$, $k=0,1,\ldots, n-1$.
For technical reasons we use estimator 
$\hat{\rho}_k$ based on $n-k$ observations and not on $n$. 
\par 
%
Let $h>0$, and for every $x\in [0, T-\delta]$ define the segment 
\[
 D_x:=\left\{ \begin{array}{ll}
               \;[x-h,x+h], &  h< x \leq T-\delta-h,\\*[2mm]
\;  [0, 2h], & 0\leq x\leq h,\\*[2mm]
\; [T-\delta-2h, T-\delta], & T-\delta -h< x\leq T-\delta.
              \end{array}
\right.
\]
Let $M_{D_x}$ be the set of indexes $k\in \{1,\ldots, n\}$ such that $k\delta\in D_x$,  
$M_{D_x}:=\{k: k\delta\in D_x\}$, and let $N_{D_x}$ be the cardinality of this set, 
$N_{D_x}:=\#\{M_{D_x}\}$.
\par 
Fix positive integer $\ell$, and 
assume that
\begin{equation}\label{eq:M-D}
 h \geq \tfrac{1}{2}(\ell + 2)\delta.
\end{equation}
For $x\in [0, T-\delta]$ let $\{a_k(x), k\in M_{D_x}\}$ denote the weights obtained as solution to the following
optimization problem
\begin{quote}
\begin{align}
 \min & \;\;\;\sum_{k\in M_{D_x}} a^2_k(x)
\nonumber
\\
\hbox{subject to } &\;\;\; \sum_{k\in M_{D_x}} a_k(x) =0,
\tag{$\sP_x$}
\\ 
&\;\;\; \sum_{k\in M_{D_x}} a_k(x) (k\delta)^j = j x^{j-1},\;\;\;\;\;j=1,\ldots, \ell.
\nonumber
\end{align}
\end{quote}
We use the convention that if $x=0$ and $j=1$ then
the right hand side of the last constraint in $(\sP_x)$ equals~$1$. 
\par 
By definition, if (\ref{eq:M-D}) holds then the linear filter associated with the weights
$\{a_k(x), k\in M_{D_x}\}$  has the following property: it reproduces without error
the first derivative of any polynomial  $p$ of ${\rm deg}(p) \leq \ell$ at point~$x$,
\begin{equation}\label{eq:weight-property}
 \sum_{k\in M_{D_x}} a_k(x) p(k\delta) = p^\prime (x),\;\;\;\forall p: \;{\rm deg}(p)\leq \ell.
\end{equation}
\par 
Now we are in a position to define our estimator of $G(x_0)$: it is given by the formula
\begin{eqnarray}
\label{eq:estimate}
\hat{G}_h(x_0) = 1 + \tfrac{1}{\lambda} \sum_{k\in M_{D_{x_0}}} a_k(x_0) \hat{R}_k,
\end{eqnarray}
where $\hat{R}_k=\hat{R}(k\delta)$, $k=0,\ldots,n-1$ are defined in (\ref{eq:R-estimator}).
\par 
The expression under the summation sign on the right hand side 
of (\ref{eq:estimate})
can be viewed as
a local polynomial estimator of the 
derivative $R^\prime(x_0)$ when the empirical covariances 
$\hat{R}_k$ are regarded as noisy observations of $R_k=R(k\delta)$. 
We refer to \citeasnoun{Gold-Nem}
for similar construction of the local polynomial estimators of derivatives
in the context of the nonparametric regression model.
\par 
The estimator $\hat{G}_h(x_0)$
depends on two design parameters, 
the {\em window width} $h$ and the 
{\em degree of polynomial} $\ell$; these parameters are specified in the sequel.
\subsection{Upper bound on the maximal risk}
Our current goal is to study accuracy of $\hat{G}_h(x_0)$. For this purpose, we 
introduce the functional class of distributions $G$ over which 
accuracy of estimator $\hat{G}_h(x_0)$ is assessed. 
\begin{definition}\label{def:smoothness-class}
\mbox{}
\begin{itemize}
\item[{\rm (i)}] Let $\beta>0$, $L>0$ be real numbers, and let 
$I\subset (0,\infty)$ be a closed interval  such that 
$x_0\in I$. 
We define $\sH_\beta (L,I)$ to be the class of all distribution functions
$G$ on $\bR_+$ such that
$G$ is  $\lfloor \beta\rfloor$ times continuously differentiable on $I$, 
and 
\[
 |G^{(\lfloor \beta\rfloor)} (x)- G^{(\lfloor \beta\rfloor)} (y)| 
\leq 
L|x-y|^{\beta-\lfloor\beta\rfloor},\;\;\;\forall 
x, y\in I;
\]
here $\lfloor \beta \rfloor$ stands for the maximal integer number strictly less than $\beta$.   
\item[{\rm (ii)}] 
We say that distribution function $G$ on $\bR_+$ belongs to the class $\sM_p(K)$, $p\geq 1$, $K>0$
if 
\[
 \rE_G [\sigma^p] = \int_0^\infty px^{p-1}[1-G(x)]\rd x \leq K <\infty.
\]
\item[{\rm (iii)}] 
Finally, we put 
\[
 \sG_\beta (L,I, K) := \sH_\beta(L, I) \cap \sM_2(K).
\]
\end{itemize}
\end{definition}
\begin{remark}\mbox{}
\begin{enumerate}
 \item[{\rm (i)}] The class $\sG_\beta(L,I,K)$ 
 imposes restrictions on smoothness
 in vicinity of $x_0$.
In all what follows the point $x_0$ 
is assumed to be fixed. If $x_0$ is separated away from zero then 
we always consider a symmetric interval $I$ centered at $x_0$: 
$I=[x_0-d,x_0+d]$ for some $0<d<x_0$. In the case $x_0=0$ we set $I=[0,2d]$. 
\item[{\rm (ii)}] The definition of $\sG_{\beta}(L,I,K)$ requires boundedness 
of the  second
moment of the service time distribution. 
This condition implies that the correlation sequence $\{H(k\delta), k\in \bZ\}$
is summable,
 which corresponds to the {\em short--term dependence} between the  values of the 
sampled discrete--time queue--length
process. 
This assumption can be relaxed. 
However, we do not pursue the case of
the long--term dependence in this paper.   
\end{enumerate}
\end{remark}
\par 
Now we are in a position to state an upper bound on the maximal risk of our estimator.
\begin{theorem}\label{th:upper-bound}
Let $x_0$ be fixed, $I:=[x_0-d, x_0+d]\subset [0, (1-\kappa )T]$ for some 
$\kappa \in (0,1)$, and suppose that $G\in \sG_\beta (L,I,K)$. 
Let $\hat{G}_*(x_0)$ be the estimator defined in 
(\ref{eq:estimate}) and associated with the degree
$\ell\geq \lfloor \beta\rfloor +1$ and the window width
\begin{equation}\label{eq:d-*}
h= h_* := \Big[
\frac{K (\sqrt{K}\vee 1)}{L^2\kappa T}\big(1+\tfrac{1}{\lambda}\big)\Big]^{1/(2\beta+2)}~.
\end{equation}
If
\begin{eqnarray}\label{eq:x-0-condition}
\frac{K(\sqrt{K}\vee 1)}{L^{2}\kappa}\big(1+\tfrac{1}{\lambda}\big)d^{-2\beta-2}
\;\leq \;T\; \leq\; 
\frac{K(\sqrt{K}\vee 1)}{L^{2}\kappa}\big(1+\tfrac{1}{\lambda}\big)\Big[\frac{2}{(\ell+2)\delta}\Big]^{2\beta+2}
\;\;\;
\end{eqnarray}
then one has
\begin{eqnarray}\label{eq:upper-bound}
\cR_{x_0}[\hat{G}_*; \sG_\beta (L,I,K) ] \;\leq\; 
C  L^{1/(\beta+1)}\Big[\frac{K (\sqrt{K}\vee 1)}{\kappa T} \big(1+\tfrac{1}{\lambda}\big)\Big]^{\beta/(2\beta+2)},  
\end{eqnarray}
where $C=C(\ell)$ depends on $\ell$ only.
\end{theorem}
\begin{remark}\label{rem:upper-bound}
\mbox{}
\begin{itemize}
 \item[{\rm (i)}] The upper bound in (\ref{eq:x-0-condition})
originates in
the requirement that the segment $D_{x_0}$ contains at least $\ell+1$ grid points.
This inequality is fulfilled if sampling is fast enough, $\delta\leq O((\kappa T)^{-1/(2\beta+2)})$.
Thus, if the asymptotics as $T\to\infty$ is considered then $\delta$ should tend to zero so that 
(\ref{eq:x-0-condition}) is fulfilled. 
The lower bound in (\ref{eq:x-0-condition}) ensures that
$D_{x_0}\subseteq I$.
\item[{\rm (ii)}] The bound in (\ref{eq:upper-bound}) is non--uniform in $x_0$; 
it is established for 
 fixed $x_0\leq (1-\kappa) T$. The bound increases as $\kappa$ gets closer to $0$ ($x_0$ approaches $T$).
This is not surprising: the empirical covariance 
estimator is not accurate  for large lags. However  
if  $x_0$ is large in comparison with $T$ then it 
is advantageous to use the trivial estimator 
$\tilde{G}(x_0)=1$. The risk of $\tilde{G}(x_0)$ admits
the following upper bound:
\begin{equation}\label{eq:upper-bound-2}
 \cR_{x_0}[\tilde{G}; \sG_\beta (L,I, K) ]  \;\leq\; Kx_0^{-2},\;\;\;\;\;\forall x_0\in \bR_+.
\end{equation}
Indeed, it follows from $G\in \sM_2(K)$ that for any $x$
\[
  1 - G(x) =\int_x^\infty \rd G(t) \leq 
x^{-2}\int_x^\infty t^2\rd G(t)
\leq Kx^{-2}.
\]
Thus, $G(x)\geq 1 - Kx^{-2}$, which  implies  (\ref{eq:upper-bound-2}).
Comparing (\ref{eq:upper-bound}) and (\ref{eq:upper-bound-2}) we see that for $x_0\leq O(T^{\beta/(4\beta+4)})$
it is advantageous to use the estimator $\hat{G}_*(x_0)$; otherwise $\tilde{G}(x_0)$ is better. 
If more stringent
conditions on the tail of $G$ are imposed [e.g., $G\in \sM_p(K)$ with $p>2$] then
the zone where $\hat{G}_*(x_0)$ is preferable becomes smaller. 
\end{itemize}
\end{remark}
\section{Estimation of arrival rate}\label{sec:arrival-rate}
The construction of Section~\ref{sec:construction} 
that led to $\hat{G}_h(x_0)$ can be used in order  
to estimate the arrival rate $\lambda$ from 
discrete observations of the queue--length process.
\par 
Let $I=[0, 2d]$ and assume that $G\in \sG_\beta(L,I,K)$. Under this condition we can 
use relation
(\ref{eq:G-R-relation}) in order to construct an estimator of $\lambda$.
Indeed, setting $t=0$ in (\ref{eq:G-R-relation}) and taking into account that $G(0)=0$
we obtain $\lambda=-R^\prime (0)$, where $R^\prime (0)$ is understood here as the 
right--side derivative of $R$ at zero. Therefore we define the estimator 
for $\lambda$
by
\begin{equation}\label{eq:lambda-estimator}
 \hat{\lambda}= - \sum_{k\in M_{D_0}} a_k(0) \hat{R}_k,
\end{equation}
where $D_0:=[0, 2h]$, $\{a_k(0), k\in M_{D_0}\}$ is the solution to ($\sP_0$)  [i.e., ($\sP_x$) with $x=0$], 
and $\hat{R}_k$, $k\in M_{D_0}$
are defined in~(\ref{eq:R-estimator}).
\par 
The next statement provides an upper bound on the risk of $\hat{\lambda}$.

\begin{theorem}\label{th:estimation-lambda}
Let $I=[0, 2d]$ and suppose that  $G\in \sG_\beta (L,I,K)$. 
Let $\hat{\lambda}_*$ denote the estimator defined in 
(\ref{eq:lambda-estimator}) and associated with degree
$\ell\geq \lfloor \beta\rfloor +1$ and window width
\begin{equation}\label{eq:d-*-arrival-rate}
h= h_* := \Big[
\frac{K (\sqrt{K}\vee 1)}{L^2 T}\Big]^{1/(2\beta+2)}~.
\end{equation}
If
\begin{eqnarray}\label{eq:x-0-condition-arrival-rate}
K(\sqrt{K}\vee 1)L^{-2}d^{-2\beta-2}\;\leq\;T\; \leq\; 
K(\sqrt{K}\vee 1)L^{-2}\Big[\frac{2}{(\ell+2)\delta}\Big]^{2\beta+2}
\end{eqnarray}
then one has
\begin{eqnarray}\label{eq:upper-bound-arrival-rate}
\sup_{G\in \sG_\beta(L,I,K)}\Big[\rE_{G,\lambda} \big| \hat{\lambda}_*-\lambda\big|^2
\Big]^{1/2} 
\;\leq\; 
C  L^{1/(\beta+1)}
 (\lambda^2+\lambda)^{1/2}\Big[\frac{K (\sqrt{K}\vee 1)}{T}\Big]^{\beta/(2\beta+2)},  
\end{eqnarray}
where $C=C(\ell)$ depends on $\ell$ only.
\end{theorem}
%
\begin{remark}\mbox{}
\begin{itemize}
\item[{\rm (i)}] The meaning of condition (\ref{eq:x-0-condition-arrival-rate}) is similar to that of
(\ref{eq:x-0-condition}), see Remark~\ref{rem:upper-bound}(i).
\item[{\rm (ii)}]
 If sampling interval $\delta$ is very small then 
 one can build an estimator which is better than 
$\hat{\lambda}_*$. 
In particular, if the  
continuous--time observation $\{X(t), 0\leq t\leq T\}$ is available
then alternative estimators of $\lambda$ can be constructed as follows
\[
 \hat{\lambda}^\uparrow =\tfrac{1}{T} \#\{t\in(0,T]: X(t)-X(t-)=1\},\;\;\;\;
\hat{\lambda}^\downarrow =\tfrac{1}{T} \#\{t\in(0,T]: X(t)-X(t-)=-1\}.
\]
Because arrivals and departures constitute  the Poisson process with intensity $\lambda$,
the mean squared errors  of $\hat{\lambda}^\uparrow$ and $\hat{\lambda}^\downarrow$ 
are given by
\[
 \rE_{G,\lambda} \big|\hat{\lambda}^\uparrow -\lambda\big|^2=
\rE_{G,\lambda} \big|\hat{\lambda}^\downarrow -\lambda\big|^2=
\lambda T^{-1},\;\;\;\forall \lambda, \forall G.
\]
Thus, in terms of dependence on the observation horizon $T$, the 
risks of  $\hat{\lambda}^\uparrow$
and $\hat{\lambda}^\downarrow$ tend to zero at the parametric rate $O(1/T)$. 
This rate is faster than the one in (\ref{eq:upper-bound-arrival-rate}).
\end{itemize}
\end{remark}
\section{Estimation of covariance  function derivative}
\label{sec:covariance}
Theorem~\ref{th:upper-bound} indicates  that 
under suitable relation between observation horizon $T$
and sampling interval $\delta$
the service time distribution $G$ can be 
estimated with the risk of the order $T^{-\beta/(2\beta+2)}$. 
In particular, for our estimator $\hat{G}_*(x_0)$ 
\[
 \cR_{x_0}[\hat{G}_*; \sG_\beta(L,I,K)]\asymp O(T^{-\beta/(2\beta+2)}),\;\;\;T\to\infty,
\]
provided that (\ref{eq:x-0-condition}) holds.
A natural question
is if this rate of convergence is optimal in the minimax sense. This is the question 
about lower bounds on the 
minimax risk $\cR^*_{x_0}[\sG_\beta(L,I,K)]$. 
\par 
Although statement~(iii) of Proposition~\ref{lem:queue-size}
provides complete probabilistic characterization of finite dimensional
distributions of the queue--length process $\{X(t), t\in\bR\}$, 
there is no explicit formula available for 
the distribution
of $X^n$. 
Because all existing techniques
for derivation of lower bounds on minimax risks rely upon sensitivity analysis of  
 the family of target distributions, such a  derivation
in the $M/G/\infty$ problem seems to be intractable.
However, some understanding of accuracy 
limitations in estimating service time distribution 
can be gained from consideration of a Gaussian approximating model. 
\par 
Proposition~\ref{lem:Gaussian-approximation} shows that 
if the arrival rate $\lambda$ is large, the
finite dimensional distributions of $\{X(t), 0\leq t\leq T\}$ are close to Gaussian.
Thus for large arrival rates  we can regard the queue--length process
as a stationary Gaussian process.  
Furthermore, 
equation (\ref{eq:G-R-relation}) shows that 
the service time distribution $G$  is proportional to the  derivative of the covariance
function of the queue--length process. 
This characterization 
suggests that,  for large arrival rates,
estimating $G$ is as hard as estimating  
derivative of the covariance function 
of a continuous--time stationary Gaussian process
from discrete observations.    
Although there is no a formal proof for statistical equivalence 
of these experiments, the assumption seems plausible.
Therefore we study the problem of estimating derivatives of covariance function
of a stationary Gaussian process from discrete observations.
\subsection{Problem formulation}
Let $X(t)$, $t\in \bR$ be a stationary Gaussian process
with zero mean and 
covariance function $\gamma\in \bL_1(\bR)$. 
The  corresponding spectral density $f$ is given by   
\begin{eqnarray*}
 f(\omega) =\int_{-\infty}^\infty \gamma(t) e^{i\omega t}\rd t= 
2\int_{0}^\infty \gamma(t) \cos(\omega t) \rd t,\;\;\;
\omega\in \bR,
\end{eqnarray*}
and, by the inverse Fourier transform,
\begin{eqnarray*}
\gamma(t) = \tfrac{1}{2\pi}\int_{-\infty}^\infty f(\omega) e^{-i\omega t}\rd \omega =\tfrac{1}{\pi}\int_0^\infty
f(\omega)\cos(\omega t)\rd \omega,\;\;\;t\in \bR.
\end{eqnarray*}
Suppose that we observe 
process $\{X(t), t\in \bR\}$ on the time interval $[0,T]$ at the points of the regular grid
$t_i=i\delta$, $i=1,\ldots, n$, where $\delta>0$ is the sampling interval, and $T=n\delta$. 
Our goal is
to estimate the first derivative, $\theta=\theta(\gamma):=\gamma^{\prime}(x_0)$, 
of $\gamma$ at fixed point $x_0\in (0,\infty)$ using the 
observation~$X^n=\{X(k\delta), k=1,\ldots, n\}$.   
\par 
Since the distribution of $X^n$ is completely determined by 
the covariance function $\gamma$ (or spectral density~$f$), 
we  write $\rP_\gamma$ and $\rE_\gamma$ for the probability measure 
and the expectation with respect to the 
distribution of $X^n$ with covariance $\gamma$. 
\par 
We  measure accuracy in estimating $\theta(\gamma)=\gamma^{\prime}(x_0)$
by the maximal risk: for any estimator $\hat{\theta}=\hat{\theta}(X^n)$ we let
\[
 \cR_{x_0}[\hat{\theta}; \sC] = \sup_{\gamma\in \sC} 
\big[\rE_\gamma \big|\hat{\theta}-\gamma^{\prime}(x_0)\big|^2\big]^{1/2},
\]
where $\sC$ is a class of target covariance functions. The minimax risk is defined by
$
\cR^*_{x_0}[\sC]=\inf_{\hat{\theta}} \cR_{x_0}[\hat{\theta};\sC]
$,
where $\inf$
is taken over all possible estimators.
\par 
In  order to relate the $M/G/\infty$ estimation problem to the present setting
let us point out
some properties of covariance functions $R(t)=\rho H(t)$  corresponding to the 
service time distributions 
$G\in \sG_\beta (L,I,K)$. First, (\ref{eq:H}) implies that if
$G\in \sH_\beta(L,I)$ [see Definition~\ref{def:smoothness-class}(i)] 
then $R\in \sH_{\beta+1}(\lambda L,I)$. 
Second, the employed moment condition $G\in \sM_2(K)$ in the $M/G/\infty$ problem boils down to 
summability of the covariance sequence
$\{R_k\}$. In the context of estimating derivative of the covariance function
this will be assumed directly. 
\par  
The above remarks  motivate the next definition. 
\begin{definition}
Let $x_0$ be fixed, and $I:=[x_0-d,x_0+d]\subset (0,\infty)$.
For $L>0$, $\beta>0$ we say that 
a covariance function $\gamma\in \bL_1(\bR)$ belongs to the functional class
$\sC_\beta(L, I, K)$ if 
\begin{itemize}
 \item[{\rm (i)}] $\int_{-\infty}^\infty |\gamma(t)|\rd t \leq K<\infty$; 
\item[{\rm (ii)}] $\gamma$ is  $\ell:=\max\{k\in \bN: k<\beta+1\}$ times continuously differentiable 
on $I$ and 
\[
 |\gamma^{(\ell)}(x)-\gamma^{(\ell)}(x^\prime)|\leq L |x-x^\prime|^{\beta+1-\ell},
\;\;\;\;\forall x, x^\prime\in I.
\]  
\end{itemize}
\end{definition}
\par 
Similarly to the definition of $\sG_\beta (L,I,K)$ in the $M/G/\infty$ estimation problem, we assume 
local smoothness around the point $x_0$ only. 
Note also that the regularity index 
of $\gamma\in \sC_\beta(L,I,K)$ equals $\beta+1$. 
We are mainly interested in bounds on the 
minimax risk $\cR^*_{x_0}[\sC_\beta(L,I, K)]$.
\subsection{Estimator and bounds on the minimax risk}
An estimator of $\theta=\theta(\gamma)=\gamma^\prime(x_0)$ 
can constructed exactly in the same way as
the estimator of $G$ in the $M/G/\infty$ problem.
Specifically, if $D_{x_0}=[x_0-h, x_0+h]$, and if  
$\{a_k(x_0), k\in M_{D_{x_0}}\}$  is the solution to optimization 
problem ($\sP_{x_0}$) then we let
\begin{equation}\label{eq:theta-h}
\hat{\theta}_h= \sum_{k\in M_{D_{x_0}}} a_k(x_0) \hat{R}_k,
\end{equation}
where $\hat{R}_k=\tfrac{1}{n-k}\sum_{i=1}^{n-k} X_i X_{i+k}$ [cf. (\ref{eq:R-estimator})].
Note that there is no need here to estimate the mean of $X(t)$ as it is assumed to be zero.
\par   
Accuracy properties of $\hat{\theta}_h$ are very similar to those of $\hat{G}_h(x_0)$.
In particular, using basically the same arguments as in the proof of Theorem~\ref{th:upper-bound}
we can establish the following result. 
\begin{theorem}\label{th:upper-theta}
 Let $I=[x_0-d, x_0+d]\subset (0,(1-\kappa)T]$ 
for some $\kappa\in (0,1)$,
and let $\gamma\in \sC_\beta (L,I,K)$. Let $\hat{\theta}_*=\hat{\theta}_{h_*}$ be the estimator
(\ref{eq:theta-h})
associated with $\ell\geq \lfloor \beta\rfloor +1$ and 
$h=h_*:= [K/(L^2\kappa T)]^{1/(2\beta+2)}$. If
\[
 KL^{-2}\kappa^{-1} d^{-2\beta-2} \;\leq\; T \;\leq\; KL^{-2}\kappa^{-1} 
\Big[\frac{2}{(\ell+2)\delta}\Big]^{2\beta+2}
\]
then 
\[
 \cR_{x_0}[\hat{\theta}_*; \sC_\beta(L,I,K)] \leq C(\ell) L^{1/(\beta+1)}
\Big(\frac{K}{\kappa T}\Big)^{\beta/(2\beta+2)}.
\]
\end{theorem}
The proof of the theorem is omitted.
\par 
Thus, the maximal risk of $\hat{\theta}_*$ 
converges to zero at the same rate as 
the risk of $\hat{G}_{h_*}(x_0)$ in the $M/G/\infty$ estimation problem;
cf. Theorem~\ref{th:upper-bound}.
\par
 The next theorem shows that this rate of convergence is, in a sense, best possible.
\begin{theorem}\label{th:lower-bound}
 Let $I=[x_0-d, x_0+d]\subset (0,\infty)$. 
There exist constants $C_1$ and $C_2$ depending on $\beta$, $x_0$, $d$ and $K$ only such that if 
\begin{equation}\label{eq:T-range}
 C_1\delta^{-2}\leq T,\;\;\;\;\;\;\;L^2T\leq C_2 \delta^{-2\beta-2}
\end{equation}
then 
\[
 \liminf_{T\to\infty} \Big\{L^{-1/(\beta+1)} T^{\beta/(2\beta+2)}\;
\cR_{x_0}^*[\sC_{\beta}(L,I, K)]\Big\} \;\geq\; c >0,
\]
where $c=c(\beta, x_0, d, K)$.
\end{theorem}
\par 
It is worth noting that the lower bound is established under condition 
$T\geq C_1\delta^{-2}$ whereas Theorems~\ref{th:upper-bound}
and~\ref{th:upper-theta} do not require it.  
We were not able to relax this condition in Theorem~\ref{th:lower-bound}.
\par 
Comparing the results of Theorems~\ref{th:upper-theta} and~\ref{th:lower-bound}
we conclude that the estimator $\hat{\theta}_*$ is rate optimal for the indicated
range of $T$ and $\delta$.  Due to relationship to the $M/G/\infty$
estimation problem, this strongly suggests that the estimator 
of the service time distribution of Section~\ref{sec:estimator} is also rate optimal. 

\section{Proofs}\label{sec:proofs}
\subsection{Proof of Proposition~\ref{lem:queue-size}}
For any  $m>1$ we write
\begin{eqnarray}\label{eq:EE}
 \rE_G \exp\Big\{\sum_{i=1}^m \theta_i X_i\Big\}= 
\rE_G \bigg\{\rE_G\Big[ \exp\Big\{\sum_{i=1}^m \theta_i X_i\Big\}
\,\Big| \{\tau_j, j\in \bZ\}\Big]\bigg\}. 
\end{eqnarray}
By (\ref{eq:X}) and by independence of $\{\tau_j, j\in \bZ\}$ and $\{\sigma_j, j\in \bZ\}$,
the conditional expectation in (\ref{eq:EE}) takes the form 
\begin{eqnarray}
 \rE_G\Big[ \exp\Big\{\sum_{i=1}^m \theta_i X_i\Big\}
\,\Big| \{\tau_j, j\in \bZ\}\Big] \;=\; \rE_G\Big[ \exp\Big\{ 
\sum_{j\in\bZ} \sum_{i=1}^m \theta_i 
{\bf 1}(\tau_j\leq t_i, \sigma_j>t_i-\tau_j)\Big\} 
\,\Big| \{\tau_j, j\in \bZ\}\Big]
\nonumber
\\
\;=\; \prod_{j\in\bZ}\rE_G\Big[ \exp\Big\{ 
\sum_{i=1}^m \theta_i {\bf 1}(\tau_j\leq t_i, \sigma_j>t_i-\tau_j)\Big\} 
\,\Big| \{\tau_j, j\in \bZ\}\Big].
\label{eq:conditional}
\end{eqnarray}
Given $x\in \bR$ consider partition of the real line by the intervals
$I_0(x)=(-\infty, t_1-x]$,
$I_k(x)=(t_k-x,t_{k+1}-x]$, $k=1,\ldots,m-1$, and $I_m(x)=(t_m-x, \infty)$.
With this notation
\begin{eqnarray*}
&& \rE_G\Big[ \exp\Big\{ 
\sum_{i=1}^m \theta_i {\bf 1}(\tau_j\leq t_i, \sigma_j>t_i-\tau_j)\Big\} 
\,\Big| \{\tau_j, j\in \rZ\}\Big] 
\\
&&\;\;\;=\;
\rP_G\{\sigma_j\in I_0(\tau_j)\} + \sum_{k=1}^m 
\exp\Big\{\sum_{i=1}^k \theta_i {\bf 1}(\tau_j\leq t_i)\Big\}
\rP_G\{\sigma_j\in I_k(\tau_j)\} 
\\
&&\;\;\;=\; 1 + \sum_{k=1}^m 
\Big[\exp\Big\{\sum_{i=1}^k \theta_i {\bf 1}(\tau_j\leq t_i)\Big\}-1\Big]
\rP_G\{\sigma_j\in I_k(\tau_j)\}. 
\end{eqnarray*}
If we let 
\[
 f(x) = \log \bigg(1+ \sum_{k=1}^m \Big[\exp\Big\{\sum_{i=1}^k \theta_i {\bf 1}(x \leq t_i)\Big\}-1\Big]
\rP_G\{\sigma_j\in I_k(x)\}\bigg),
\]
then in view of (\ref{eq:EE}), (\ref{eq:conditional}) and 
Campbell's theorem \cite[Section~3.2]{Kingman} we obtain 
\[
\rE_G \exp\Big\{\sum_{i=1}^m \theta_i X_i\Big\}=
\rE_G \exp\Big\{\sum_{j\in \bZ} f(\tau_j)\Big\}=
\exp\bigg\{\lambda \int_{-\infty}^\infty [e^{f(x)} -1]\rd x\bigg\}.
\]
\par 
Denote $S_m(\theta)=\mu \int_{-\infty}^\infty [e^{f(x)}-1]\rd x$;  
our current goal is to compute this 
integral. 
We have
\begin{eqnarray}
\int_{-\infty}^\infty [e^{f(x)}-1]\rd x 
&=&  \sum_{k=1}^{m-1} \int_{-\infty}^\infty 
 \big(\exp\big\{\sum_{i=1}^k \theta_i {\bf 1}(x \leq t_i)\big\}-1\big)
\big[\bar{G}(t_k-x)- \bar{G}(t_{k+1}-x)\big]\rd x
\nonumber
\\
&&\;\;+\;\; \int_{-\infty}^\infty \big(\exp\big\{\sum_{i=1}^m
\theta_i {\bf 1}(x \leq t_i)\big\}-1\Big)
\bar{G}(t_m-x)\rd x
\nonumber
\\
&
=:&\;\sum_{k=1}^{m-1} J_k + L_m,
\nonumber
\end{eqnarray}
where we denoted for brevity $\bar{G}=1-G$.
For $k=1,\ldots, m-1$ we obtain
\begin{eqnarray}
&&J_k = \big( \exp\big\{\sum_{i=1}^k \theta_i\big\} -1 \big) 
\int_{-\infty}^{t_1}
\big[\bar{G}(t_k-x)- \bar{G}(t_{k+1}-x)\big]\rd x
\nonumber
\\
&&\;\;\;\;\;\;\;\;\;+\;
\sum_{j=1}^{k-1}
\big(\exp\big\{\sum_{i=j+1}^k \theta_i\big\} - 1\big)
\int_{t_j}^{t_{j+1}} \big[\bar{G}(t_k-x)- \bar{G}(t_{k+1}-x)\big] \rd x
\nonumber
\\
&&\;\;\;\;\;\;=\; \tfrac{1}{\mu}\big( \exp\big\{\sum_{i=1}^k \theta_i\big\} -1 \big)
\big[H(t_k-t_1)-H(t_{k+1}-t_1)\big]
\nonumber
\\
&&\;\;\;\;\;\;\;\;\;+\;\tfrac{1}{\mu} \sum_{j=1}^{k-1}
\big(\exp\big\{\sum_{i=j+1}^k \theta_i\big\} - 1\big) \big[ H(t_k-t_{j+1})-H(t_k-t_j)-H(t_{k+1}-t_{j+1})
+ H(t_{k+1}-t_j)\big]
\nonumber
\\
&&\;\;\;\;\;\;=\;\tfrac{1}{\mu}\big( \exp\big\{\sum_{i=1}^k \theta_i\big\} -1 \big)
\big[H_{k-1}-H_{k}\big]
\nonumber
\\
&&\;\;\;\;\;\;\;\;\;\;\;\;+\;\;
\tfrac{1}{\mu} \sum_{j=1}^{k-1}
\big(\exp\big\{\sum_{i=j+1}^k \theta_i\big\} - 1\big) 
[ H_{k-j-1}-2H_{k-j}+ H_{k-j+1}].
\label{eq:S-k}
\end{eqnarray}
Similarly, 
\begin{eqnarray}
L_m &=& 
\tfrac{1}{\mu}\big( \exp\big\{\sum_{i=1}^m \theta_i\big\} -1 \big) H(t_m-t_1)
\;+\; \tfrac{1}{\mu} \sum_{j=1}^{m-1}
\big(\exp\big\{\sum_{i=j+1}^m \theta_i\big\} - 1\big) 
\big[ H(t_m-t_{j+1})-H(t_m-t_j)\big]
\nonumber
\\
&=&
\tfrac{1}{\mu}\big( \exp\big\{\sum_{i=1}^m \theta_i\big\} -1 \big) H_{m-1}
+
\tfrac{1}{\mu} \sum_{j=1}^{m-1}
\big(\exp\big\{\sum_{i=j+1}^m \theta_i\big\} - 1\big) 
\big[ H_{m-j-1}-H_{m-j}\big].
\label{eq:S-d}
\end{eqnarray}
The usual convention $\sum_{k=j}^m =0$ 
if $m<j$ is employed in  (\ref{eq:S-k}) and (\ref{eq:S-d})
and from now on.
\par
Note that by definition $S_m(\theta)=\mu \sum_{k=1}^{m-1}J_k +\mu L_m$, and 
we have the following recursive formula 
\begin{eqnarray}\label{eq:S-iteration}
 S_{m+1}(\theta)=S_{m}(\theta)+ \mu (J_m - L_m + L_{m+1}).
\end{eqnarray}
For any $m>1$, using (\ref{eq:S-k}) and (\ref{eq:S-d}),
after straightforward algebraic manipulations we obtain 
\begin{eqnarray*}
&& \mu (J_m - L_m + L_{m+1})
 \\
&& \;= \;\big(
 e^{\sum_{i=1}^n \theta_i}-1\big)(H_{m-1}-H_m)
\;+\;\sum_{j=1}^{m-1} 
\big(e^{\sum_{i=j+1}^m \theta_i}-1\big)(H_{m-j-1}-2H_{m-j}+H_{m-j+1})
\\
&&\;\;\; -\; \big(e^{\sum_{i=1}^m\theta_i}-1\big)H_{m-1} \;-\; \sum_{j=1}^{m-1}
\big(e^{\sum_{i=j+1}^m \theta_i} -1\big)(H_{m-j-1}-H_{m-j})
\\
&&\;\;\;+\; \big(e^{\sum_{i=1}^{m+1}\theta_i}-1\big)H_m \;+\;
\sum_{j=1}^m \big(e^{\sum_{i=j+1}^{m+1}\theta_i}-1\big)(H_{m-j}-H_{m-j+1})
\\
&&\;=\;(e^{\theta_{m+1}}-1) \;+\;
\sum_{k=1}^m H_k (e^{\theta_{m-k+1}}-1) 
e^{\sum_{i=m-k+2}^m \theta_i} 
(e^{\theta_{m+1}}-1).
 \end{eqnarray*}
Taking into account that $S_1(\theta)=e^{\theta_1}-1$ and
iterating the formula (\ref{eq:S-iteration})
we obtain
\begin{eqnarray*}
 S_{n+1}(\theta) &=& (e^{\theta_1}-1) + \sum_{m=1}^{n} (e^{\theta_{m+1}}-1)
 \;+\; 
 \sum_{m=1}^n \sum_{k=1}^m
 H_k(e^{\theta_{m-k+1}}-1)e^{\sum_{i=m-k+2}^m\theta_i}(e^{\theta_{m+1}}-1)
 \\
 &=& 
 \sum_{m=1}^{n+1} (e^{\theta_m}-1) +
 \sum_{k=1}^n H_k \sum_{m=k}^n 
 (e^{\theta_{m-k+1}}-1)e^{\sum_{i=m-k+2}^m\theta_i}(e^{\theta_{m+1}}-1).
\end{eqnarray*}
This completes the proof. 
\epr
\subsection{Proof of Proposition~\ref{cor:moments}}
The proof involves straightforward though tedious differentiation of 
(\ref{eq:4-var}). 
\par 
Let $S(\theta)$ stand for 
the right hand side of (\ref{eq:4-var}), where $(\theta_1,\theta_2,\theta_3,\theta_4)$ 
is replaced by $(\theta_i, \theta_j, \theta_k, \theta_m)$ for convenience.
Denote 
$\psi(\theta)=\rE_G\exp\{\theta_i X_i +\theta_j X_j + \theta_k X_k +\theta_mX_m\}$. 
It is checked by direct calculation that 
\begin{eqnarray}\label{eq:a-a}
 \frac{\partial^4 \psi(\theta)}{\partial\theta_i\partial\theta_j
\partial\theta_k\partial\theta_m} = \exp\{-\rho S(\theta)\} \big[a_1(\theta)\rho + 
a_2(\theta)\rho^2+a_3(\theta)\rho^3+ a_4(\theta)\rho^4\big],
\end{eqnarray}
where $a_1(\theta), a_2(\theta), a_3(\theta)$ and $a_4(\theta)$ are given by the following expressions:
\begin{eqnarray*}
 a_1(\theta)&=& S_{\theta_i\theta_j\theta_k\theta_m},
\\
 a_2(\theta) &=& S_{\theta_i\theta_j\theta_k}S_{\theta_m}+S_{\theta_i\theta_j\theta_m}S_{\theta_k}
+S_{\theta_j\theta_k\theta_m}S_{\theta_i} + S_{\theta_i\theta_k\theta_m} S_{\theta_j} + S_{\theta_i\theta_j}
S_{\theta_k\theta_m} + S_{\theta_i\theta_k}S_{\theta_j\theta_m}+ S_{\theta_j\theta_k}S_{\theta_i\theta_m},
\\
 a_3(\theta)&=&S_{\theta_i\theta_j}S_{\theta_k}
S_{\theta_m} + S_{\theta_i\theta_k} S_{\theta_j}S_{\theta_m} + 
S_{\theta_i\theta_m} S_{\theta_j}S_{\theta_k}
+
S_{\theta_j\theta_k}S_{\theta_i}S_{\theta_m}
+S_{\theta_j\theta_m}S_{\theta_i}S_{\theta_k} + S_{\theta_k\theta_m} S_{\theta_i}S_{\theta_j},
\\
a_4(\theta) &=& S_{\theta_i}S_{\theta_j}S_{\theta_k}S_{\theta_m}.
\end{eqnarray*}
Here  we put for brevity $S_{\theta_{j_1}\cdots \theta_{j_k}}=S_{\theta_{j_1}\cdots \theta_{j_k}}(\theta):=
\partial^k S(\theta)/\partial \theta_{j_1}
\cdots\partial \theta_{j_k}$.  
In fact, expression (\ref{eq:a-a}) is obtained by application of di~Bruno's formula 
for derivatives of composite functions 
[see, e.g., \citeasnoun[Chapter~2]{Riordan}] to  (\ref{eq:4-var}).
\par 
In order to complete the proof, it is sufficient to note that
\begin{equation}\label{eq:S-0}
S(0)=1,\;\; 
S_{\theta_j}(0)=1,\;\; \forall j, 
\end{equation}
and  for any $j_1\leq j_2\leq j_3\leq j_4$ 
\begin{eqnarray}\label{eq:S-0-0} 
S_{\theta_{j_1}\theta_{j_2}}(0) = H_{j_2-j_1},\;
\;S_{\theta_{j_1}\theta_{j_2}\theta_{j_3}}(0)  = H_{j_3-j_1}, \;\;
S_{\theta_{j_1}\cdots\theta_{j_4}}(0)= H_{j_4-j_1}.
\end{eqnarray}
\par 
Although (\ref{eq:S-0-0}) is proved for $1\leq i\leq j\leq k\leq m\leq n$, 
a similar result holds more generally.
With the introduced definition of $q_I$,
(\ref{eq:S-0}), (\ref{eq:S-0-0}) imply that
\begin{eqnarray*}
a_4(0)&=& 1
\\
a_3(0)&=& H_{|i-j|} + H_{|k-i|} + H_{|m-i|} + H_{|k-j|} + H_{|m-j|} + H_{|m-k|}
\\
a_2(0) &=&  H_{q_{\{i,j,k\}}} + H_{q_{\{i,j,m\}}}+ H_{q_{\{j,k,m\}}} + H_{q_{\{i,k,m\}}}
+ H_{|i-j|}H_{|k-m|} + H_{|k-i|}h_{|m-j|} + H_{|k-j|}H_{|m-i|}
\\
a_1(0) &=& H_{q_{\{i,j,k,m\}}}.
\end{eqnarray*}
This completes the proof. 
\epr
\subsection{Proof of Theorems~\ref{th:upper-bound} and~\ref{th:estimation-lambda}}
Throughout the proof $C_i, c_i$, $i=1,2,\ldots$ stand for constants depending on $\ell$ only, 
unless it is mentioned
explicitly. The proofs of both theorems are almost identical. We first prove Theorem~\ref{th:upper-bound}
and then indicate modifications needed for the proof of Theorem~\ref{th:estimation-lambda}.
\par
It follows from (\ref{eq:G-R-relation}) and (\ref{eq:estimate}) that
\begin{eqnarray*}
 \hat{G}_h(x_0)- G(x_0) 
&=& \tfrac{1}{\lambda}\Big[\sum_{k\in M_{D_{x_0}}} a_k(x_0) \hat{R}_k - R^\prime(x_0)\Big].
\\
&=& \tfrac{1}{\lambda}\Big[\sum_{k\in M_{D_{x_0}}} a_k(x_0) (\hat{R}_k - R_k) +
\sum_{k\in M_{D_{x_0}}} a_k(x_0) R_k - R^\prime(x_0)\Big].
\end{eqnarray*}
Therefore 
\begin{eqnarray}
&& \Big[\rE_G  |\hat{G}_h(x_0)- G(x_0)|^2\Big]^{1/2}
\nonumber
\\
&& \;\;\;\;\;\leq\; \tfrac{1}{\lambda} 
\Big\{\rE_G \Big[\sum_{k\in M_{D_{x_0}}} a_k(x_0) (\hat{R}_k - R_k)\Big]^2\Big\}^{1/2} + 
\tfrac{1}{\lambda}\Big|\sum_{k\in M_{D_{x_0}}} a_k(x_0) R_k - R^\prime(x_0)\Big|. 
\label{eq:bias-variance}
\end{eqnarray}
In the subsequent proof
we bound the expression on the right hand side of the above display formula. 
The result of the theorem will follow
from series of lemmas given below.  
\par\medskip 
We begin with a well known result on the  properties of the local polynomial estimators; 
see, e.g., \citeasnoun[Lemma~1.3.1]{N00} and \citeasnoun[Section~1.6]{tsybakov}.
\begin{lemma}\label{lem:norm-a}
 Let 
$\{a_k(x_0),\;k\in M_{D_{x_0}}\}$ be the solution to ($\sP_x$), and let  
(\ref{eq:M-D}) hold; then
\begin{equation}\label{eq:a-coefficients}
 \Big[
\sum_{k\in M_{D_{x_0}}}
 |a_k(x_0)|^2\Big]^{1/2} \;\leq\; \frac{C_1}{h\sqrt{N_{D_{x_0}}}}, \;\;\;\;\;
\sum_{k\in M_{D_{x_0}}} |a_k(x_0)| \leq \frac{C_2}{h},
\end{equation}
where $C_1=C_1(\ell)$ and $C_2=C_2(\ell)$ are constants  depending on $\ell$ only.
\end{lemma}
\par 
The next result establishes an upper bound on accuracy of the empirical covariance estimator.
\begin{lemma}\label{lem:accuracy-of-R-k}
For any $k=0,\ldots, n-1$ one has
\[
\rE_G |\hat{R}_k - R_k|^2 \leq \tfrac{C_3}{n-k} (\rho^2+\rho) \sum_{i=1}^n H_i,
\]
where $C_3$ is an absolute constant.
\end{lemma}
\pr  Let
$\tilde{R}_k := \tfrac{1}{n-k} \sum_{i=1}^{n-k} (X_i-\rho)(X_{i+k}-\rho)$; then
$\rE_G \tilde{R}_k = R_k$, and by definition of $\hat{\rho}_k$ 
\begin{eqnarray*}
 \hat{R}_k = \tfrac{1}{n-k}\sum_{i=1}^{n-k} (X_i-\hat{\rho}_k)(X_{i+k}-\hat{\rho}_k)= 
\tilde{R}_k - (\hat{\rho}_k -\rho)^2.
\end{eqnarray*}
Therefore 
\begin{eqnarray}
\rE_G |\hat{R}_k - R_k|^2 =
\rE_G |\tilde{R}_k - R_k|^2 - 2 \rE_G \big[\tilde{R}_k(\hat{\rho}_k - \rho)^2\big] + \rE_G |\hat{\rho}_k - \rho|^4
=: J_1 - 2 J_2 + J_3.
\label{eq:cov-00}
\end{eqnarray}
Now we proceed with computation of the  terms 
on the right hand side of (\ref{eq:cov-00}). 
\par\medskip 
1$^0$.  {\em Computation of $J_1$.}
\par
Let $r_k := \rE_G [X_i X_{i+k}]=R_k+\rho^2=\rho H_k+\rho^2$ and $\tilde{r}_k:=\frac{1}{n-k} \sum_{i=1}^{n-k}X_i X_{i+k}$; 
then  
\begin{eqnarray*}
 \tilde{R}_k - R_k = \tilde{r}_k - r_k + 2\rho^2 - \tfrac{\rho}{n-k}\sum_{i=1}^{n-k}(X_i+X_{i+k}).
\end{eqnarray*}
Thus
\begin{eqnarray}
 J_1 &=& \rE_G|\tilde{R}_k-R_k|^2 
\nonumber
\\
&=& \rE_G|\tilde{r}_k - r_k|^2 -2
\rE_G \Big [ (\tilde{r}_k-r_k) \tfrac{\rho}{n-k}\sum_{i=1}^{n-k} (X_i+X_{i+k})\Big] +
\rE_G\Big[2\rho^2 -\tfrac{\rho}{n-k}\sum_{t=1}^{n-k}(X_i+X_{i+k})\Big]^2
\nonumber
\\
&=:& J_1^{(1)} - J_1^{(2)} + J_{1}^{(3)}.
\label{eq:J-1-0}
\end{eqnarray}
Equality (\ref{eq:moments-1}) of Proposition~\ref{cor:moments}  implies that for
any $k = 0,\ldots, n$ and $i, j=1,\ldots, n-k$ one has
\begin{eqnarray*}
  \rE_G \big[X_iX_{i+k}X_jX_{j+k}\big]&=& \rho^4 + 
\rho^3\big[H_k+H_{|j-i|} + H_{|j-i+k|}+ H_{|j-i-k|}+H_{|j-i|} + H_k\big]
\\*[2mm]
&&+\;\rho^2\big[2H_{k\vee |j-i|\vee |j-i-k|}+2H_{k\vee|j-i|\vee|j-i+k|} +
H_k^2+H^2_{|j-i|}+H_{|j-i+k|}H_{|j-i-k|}\big]
\\*[2mm]
&&+\;\rho H_{k\vee|j-i|\vee|j-i+k|\vee|j-i-k|}.
\end{eqnarray*}
Since
$r_k^2=\rho^4+2\rho^3H_k+\rho^2 H_k^2$,
\begin{eqnarray}
J_1^{(1)}&=&\rE_G |\tilde{r}_k-r_k|^2 
= \tfrac{1}{(n-k)^2} \sum_{i,j =1}^{n-k} \rE_G \big[X_iX_{i+k}X_jX_{j+k}\big] - r_k^2
\nonumber 
\\
&=&\tfrac{1}{(n-k)^2}\sum_{i, j=1}^{n-k}
\Big\{ \rho^3 \big[2H_{|j-i|}+ H_{|j-i-k|}+H_{|j-i+k|}\big] \;+\;\rho H_{k\vee|j-i|\vee|j-i-k|\vee|j-i+k|}
\label{eq:J-1-1}
\\
&&\hspace{15mm}\;\;\;+\;\rho^2\big[2H_{k\vee|j-i|\vee|j-i-k|} + 2H_{k\vee|j-i|\vee|j-i+k|} + 
H^2_{|j-i|} + H_{|j-i-k|}H_{|j-i+k|}\big]\Big\}.
\nonumber
\end{eqnarray}
Furthermore, 
\begin{eqnarray}
 J_1^{(3)} &=& 4\rho^2 - \tfrac{4\rho^3}{n-k} \rE_G \sum_{i=1}^{n-k} (X_i + X_{i+k})
+ \rE_G\tfrac{\rho^2}{(n-k)^2} \sum_{i,j=1}^{n-k} (X_i+X_{i+k})(X_j+X_{j+k})
\nonumber 
\\
&=& 
-4\rho^4 + \tfrac{\rho^2}{(n-k)^2} \sum_{i,j=1}^{n-k} \big[2r_{|j-i|} + r_{|j-i+k|}+r_{|j-i-k|}\big]
\nonumber
\\
&=& \tfrac{\rho^3}{(n-k)^2} \sum_{i,j=1}^{n-k}\big[2H_{|j-i|}+H_{|j-i-k|}+H_{|j-i+k|}\big].
\label{eq:J-1-3}
\end{eqnarray}
Now we proceed with $J_{1}^{(2)}$:
\begin{eqnarray*}
J_1^{(2)} &=& \tfrac{2\rho}{n-k}\sum_{i=1}^{n-k} \rE_G (\tilde{r}_k-r_k)(X_i+X_{i+k}) =
\tfrac{2\rho}{n-k} \sum_{i=1}^{n-k} \big[ \rE_G (\tilde{r}_k X_i + \tilde{r}_k X_{i+k}) - 2\rho(\rho^2+\rho H_k)\big].
\end{eqnarray*}
We have
\begin{eqnarray*}
\rE_G [\tilde{r}_k X_i] &=& \tfrac{1}{n-k}\sum_{j=1}^{n-k} \rE_G[X_j X_{j+k} X_i]
\\
&=& \tfrac{1}{n-k}\sum_{j=1}^{n-k} \big[\rho^3+\rho^2 H_k + \rho^2 H_{|j-i|} + \rho^2 H_{|j-i-k|} +
\rho H_{k\vee |j-i|\vee |i-j-k|}\big]
\\
\rE_G [\tilde{r}_k X_{i+k}] &=& \tfrac{1}{n-k}\sum_{j=1}^{n-k} \rE_G[X_j X_{j+k} X_{i+k}]
\\
&=& \tfrac{1}{n-k}\sum_{j=1}^{n-k} \big[\rho^3+\rho^2 H_k + \rho^2 H_{|j-i|} + \rho^2 H_{|j-i+k|} +
\rho H_{k\vee |j-i|\vee |i-j+k|}\big],
\end{eqnarray*}
which yields
\begin{eqnarray}
J_1^{(2)}
= \tfrac{2}{(n-k)^2} \sum_{i,j=1}^{n-k} \big[2\rho^3 H_{|j-i|}+ \rho^3 H_{|i-j+k|}+ \rho^3H_{|i-j-k|} 
+ \rho^2H_{|i-j|\vee k\vee |i-j-k|} + \rho^2H_{|i-j|\vee k\vee |i-j+k|}\big].
\label{eq:J-1-2}
\end{eqnarray}
Combining (\ref{eq:J-1-2}), (\ref{eq:J-1-3}), (\ref{eq:J-1-1}) and (\ref{eq:J-1-0}) we obtain 
\begin{eqnarray*}\label{eq:J-1}
 J_1= \tfrac{\rho^2}{(n-k)^2} \sum_{i,j=1}^{n-k} 
\big[H_{|i-j|}^2+H_{|i-j+k|}H_{|i-j-k|}\big] + \tfrac{\rho}{(n-k)^2}\sum_{i,j=1}^{n-k} 
H_{k\vee |i-j|\vee |i-j-k|\vee |i-j+k|}.
\end{eqnarray*}
Taking into account that $H$ is a monotone decreasing function, and $H(0)=1$ we obtain   
\[
 J_1 \leq \tfrac{c_1}{n-k} (\rho^2+\rho)\sum_{i=1}^n H_i,
\]
where $c_1$ is an absolute constant.

\par\medskip 
2$^0$. {\em Computation of $J_2$.}
It follows from the definition of $J_2$ that
\[
 J_2 = \rho^3 H_k - 2\rho \rE_G[\tilde{R}_k\hat{\rho}_k] + \rE_G[\tilde{R}_k\hat{\rho}_k^2].
\]
We have 
\begin{eqnarray*}
 E_G \big[\tilde{R}_k \hat{\rho}_k\big] &=& \tfrac{1}{(n-k)^2}\, \rE_G \sum_{i,j=1}^{n-k} (X_i-\rho)(X_{i+k}-\rho) X_{j}
\\
&=&
\rE_G \tfrac{1}{(n-k)^2}\,\sum_{i,j=1}^{n-k} \Big[X_iX_{i+k}X_{j} - \rho X_{i+k}X_{j} - \rho X_iX_{j} + 
\rho^2 X_{j}
\Big]
\\
&=&
\rho^2 H_k +
\tfrac{\rho}{(n-k)^2} \sum_{i,j=1}^{n-k} H_{k\vee|i-j|\vee |i-j+k|}~.
\end{eqnarray*}
Furthermore, 
\begin{eqnarray*}
&&\rE_G [\tilde{R}_k\hat{\rho}_k^2] = \tfrac{1}{(n-k)^3} \sum_{i,j,l=1}^{n-k}
\rE_G\Big[X_iX_{i+k}X_jX_l - \rho X_{i+k}X_jX_l - \rho X_iX_jX_l+\rho^2 X_jX_l\Big]
\\
&&=  \tfrac{\rho^2}{(n-k)^3} \sum_{i,j,l=1}^{n-k}
\Big[ H_{k\vee |i-j|\vee |i-j+k|} + H_{|i-l|\vee|i-l+k|\vee k} + H_k H_{|j-l|} +
H_{|i-j|}H_{|i-l+k|}+H_{|i-l|}H_{|i-j+k|}\Big]
\\
&& \hspace{60mm}
 + \;\rho^3 H_k +\tfrac{\rho}{(n-k)^3} \sum_{i,j,l=1}^{n-k} H_{k\vee|i-j|\vee|i-l|\vee|i+k-j|\vee|i+k-l|\vee|j-l|}
\\
&&= \rho^3 H_k + \tfrac{\rho^2}{(n-k)^2} \sum_{i,j=1}^{n-k} 2H_{k\vee |i-j|\vee |i-j+k|}
\\
&& + \tfrac{1}{(n-k)^3} \sum_{i,j,l=1}^{n-k}
\Big[  \rho^2\big(H_k H_{|j-l|} +
H_{|i-j|}H_{|i-l+k|}+H_{|i-l|}H_{|i-j+k|}\big) + \rho  H_{k\vee|i-j|\vee|i-l|\vee|i+k-j|\vee|i+k-l|\vee|j-l|}
\Big].
\end{eqnarray*}
Combining these equalities we obtain
\begin{eqnarray*}
 J_2 &=&\tfrac{1}{(n-k)^3} \sum_{i,j,l=1}^{n-k}
\Big[  \rho^2\big(H_k H_{|j-l|} +
H_{|i-j|}H_{|i-l+k|}+H_{|i-l|}H_{|i-j+k|}\big) + \rho  H_{k\vee|i-j|\vee|i-l|\vee|i+k-j|\vee|i+k-l|\vee|j-l|}\Big]
\\
&\leq & \tfrac{c_2}{n-k}(\rho^2+\rho) \sum_{i=1}^{n} H_i,
\end{eqnarray*}
where $c_2$ is an absolute constant.
\par\medskip 
3$^0$. {\em Computation of $J_3$.}
By definition,
$J_3= \rE_G|\tfrac{1}{n-k}\sum_{i=1}^{n-k} (X_i -\rho)|^4$. Using Proposition~\ref{cor:moments} 
after routine calculations
we obtain for all $i,j,l,m=1,\ldots, n-k$
\begin{eqnarray*}
 &&\rE_G\big[(X_i-\rho)(X_j-\rho)(X_l-\rho)(X_m-\rho)\big] 
\\
&&\;=\;\rho^2\Big[H_{|i-j|}H_{|l-m|}+ H_{|i-l|}H_{|j-m|}+H_{|l-j|}H_{|i-m|}\big] + \rho 
H_{|i-j|\vee|j-l|\vee|l-m|\vee|l-m|\vee|i-m|\vee|i-l|},
\end{eqnarray*}
so that 
\begin{eqnarray*}
 J_3&=& \tfrac{\rho^2}{(n-k)^4} \sum_{i,j,l,m=1}^{n-k}
\Big[H_{|i-j|}H_{|l-m|}+ H_{|i-l|}H_{|j-m|}+H_{|l-j|}H_{|i-m|}\big] 
\\
&&\; 
+ \; \tfrac{\rho}{(n-k)^4} \sum_{i,j,l,m=1}^{n-k}
H_{|i-j|\vee|j-l|\vee|l-m|\vee|l-m|\vee|i-m|\vee|i-l|} \;\leq\;\tfrac{c_3}{n-k} (\rho^2+\rho) \sum_{i=1}^n H_i, 
\end{eqnarray*}
where $c_3$ is an absolute constant.
\par 
Combining inequalities for $J_1$, $J_2$ and $J_3$ with   (\ref{eq:cov-00}) we complete the proof.
\epr 
\begin{lemma}\label{lem:stoch-error}
For every $x_0\in [0, T-\delta]$ one has  
\begin{eqnarray*}
\rE_G \Big|\sum_{k\in M_{D_{x_0}}} a_k(x_0) (\hat{R}_k - R_k)\Big|^2 \leq  
\frac{C_4\delta}{h^{2} \psi_{x_0}(T)}\big(\rho^2+\rho\big) \sum_{i=1}^n H_i,
\end{eqnarray*}
where $C_4=C_4(\ell)$ is a constant depending on $\ell$ only, and 
\[
 \psi_{x_0}(T)=\psi_{x_0}(T, h, \delta):=\left\{\begin{array}{ll}
                        T-x_0-h, & h \leq x_0< T-\delta-h,\\
                         T-2h, & 0 \leq x_0\leq h,\\
                         \delta, & T-\delta - h\leq x_0\leq T-\delta.  
                       \end{array}
\right.
\]
\end{lemma}
\pr 
By Lemmas~\ref{lem:norm-a} and~\ref{lem:accuracy-of-R-k} and by the 
Cauchy--Schwarz inequality 
\begin{eqnarray}\label{eq:sq-bound}
\rE \Big|\sum_{k\in M_{D_{x_0}}} a_k(x_0) (\hat{R}_k - R_k)\Big|^2 \leq 
\sum_{k\in M_{D_{x_0}}} a_k^2(x_0) 
\sum_{k\in M_{D_{x_0}}} \rE_G (\hat{R}_k-R_k)^2
\nonumber
\\
\leq 
\frac{c_1^2}{h^2N_{D_{x_0}}}\big(\rho^2+\rho\big) \Big(\sum_{i=1}^n H_i\Big)
\sum_{k\in M_{D_{x_0}}}\tfrac{1}{n-k}.
\end{eqnarray}
Let $\underline{k}=\min\{k\in (1,\ldots, n-1): k\in M_{D_{x_0}}\}$ and 
$\overline{k}=\max\{k\in (1,\ldots, n-1): k\in M_{D_{x_0}}\}$; then 
\begin{eqnarray*}
 \sum_{k\in M_{D_{x_0}}} \tfrac{1}{n-k} = 
\sum_{k=\underline{k}}^{\overline{k}} \tfrac{1}{n-k} 
\leq \ln \Big(\frac{n-\underline{k}}{n-\overline{k}}\Big) = 
\ln \Big(1+ 
\frac{\overline{k}-\underline{k}}{n-\overline{k}}\Big) \leq \frac{\overline{k}-\underline{k}}{n-\overline{k}}.
\end{eqnarray*}
\par 
First, assume that $D_{x_0}=[x_0-h, x_0+h]$. In this case
$\underline{k}=[(x_0-h)/\delta]+1$, $\overline{k}=[(x_0+h)/\delta]$, where $[\cdot]$ is the integer 
part, and then
$ \sum_{k\in M_{D_{x_0}}} 1/(n-k) \;\leq\; 2h/(T - x_0-h)$.
If $D_{x_0}=[0, 2h]$ then $\underline{k}=1$ and $\overline{k}=[2h/\delta]$ which leads to 
 $\sum_{k\in M_{D_{x_0}}} 1/(n-k) \leq 2h/(T - 2h)$.
Finally, if $D_{x_0}=[T-2h-\delta, T-\delta]$ then
$\bar{k}=n-1$, $\underline{k}=(n-1)-[2h/\delta]$, and 
$\sum_{k\in M_{D_{x_0}}} 1/(n-k) \leq 2h/\delta$.
\par 
Combining these bounds with (\ref{eq:sq-bound}) and 
taking into account that $(2h/\delta)-1\leq N_{D_{x_0}}\leq (2h/\delta)+1$, we complete the proof.
\epr 
\begin{lemma}\label{lem:bias}
 Let $G\in \sH_\beta(L, I)$, $I=[x_0-d,x_0+d]\supseteq D_{x_0}$, 
and $\{a_k(x_0), k\in M_{D_{x_0}}\}$ be the weights defined
by ($\sP_{x_0}$) with $\ell\geq \lfloor\beta\rfloor+1$. Assume that
(\ref{eq:M-D}) holds; then 
\[
\Big|\sum_{k\in M_{D_{x_0}}} a_k(x_0) R_k - R^\prime(x_0)\Big|\;\leq\; C_2\lambda L h^{\beta},
\]
where $C_2=C_2(\ell)$ is the constant appearing in (\ref{eq:a-coefficients}).
\end{lemma}
\pr 
Recall $R(t)=\rho^2+\rho h(t)=\rho^2 +\lambda \int_t^\infty [1- G(x)]\rd x$; this implies
\[
R^\prime(t)=-\lambda(1-G(t)),\;\;\;
R^{(j)}(t)=\lambda G^{(j-1)}(t),\;\; \forall j=2,\ldots, \lfloor\beta\rfloor+1.
\]
Thus if $G\in \sH_\beta(L, I)$ then $R\in \sH_{\beta+1}(\lambda L, I)$. 
Since $D_{x_0}\subseteq I$, function $R$ can be expanded in the Taylor series around $x_0$. 
In particular, for any $k\in M_{D_{x_0}}$
\begin{equation}\label{eq:R-1}
R(k\delta)=R(x_0)+\sum_{j=1}^{\lfloor\beta\rfloor} \tfrac{1}{j!} R^{(j)}(x_0)(k\delta -x_0)^j + 
\tfrac{1}{(\lfloor \beta\rfloor +1)!} R^{(\lfloor \beta\rfloor+1)}(\xi_k)(k\delta -x_0)^{\lfloor \beta\rfloor +1},
\end{equation}
where $\xi_k=\tau k\delta + (1-\tau)x_0$ for some $\tau\in [0,1]$.
Denote 
\begin{equation}\label{eq:R-2}
\bar{R}_{x_0}(y) := R(x_0) + \sum_{j=1}^{\lfloor\beta\rfloor +1} \tfrac{1}{j!} 
R^{(j)}(x_0) (y-x_0)^j,\;\;\;y\in D_{x_0}. 
\end{equation}
Because $\bar{R}_{x_0}(\cdot)$ is a polynomial of degree $\lfloor\beta\rfloor+1$
and  $\ell\geq \lfloor \beta\rfloor +1$, we have by (\ref{eq:weight-property}) that  
\[
 \sum_{k\in M_{D_{x_0}}} a_k(x_0) \bar{R}_{x_0}(k\delta)=\bar{R}_{x_0}^\prime(x_0)=R^\prime(x_0).
\]
Therefore
\begin{eqnarray*}
&& \sum_{k\in M_{D_{x_0}}} a_k(x_0) R(k\delta)  - R^\prime(x_0) \;=\; \sum_{k\in M_{D_{x_0}}}
a_k(x_0)\big[R(k\delta) - \bar{R}_{x_0}(k\delta)\big]
\\
&&\;\;\;=\;
\sum_{k\in M_{D_{x_0}}} \tfrac{1}{(\lfloor\beta\rfloor +1)!}\,a_k(x_0) 
\big[R^{(\lfloor \beta\rfloor+1)}(\xi_k) -R^{(\lfloor \beta\rfloor+1)}(x_0)\big]  
(k\delta-x_0)^{\lfloor\beta\rfloor +1},
\end{eqnarray*}
where we have used (\ref{eq:R-1}) and (\ref{eq:R-2}).
This yields 
\begin{eqnarray*}
 \Big|\sum_{k\in M_{D_{x_0}}} a_k(x_0) R(k\delta)  - R^\prime(x_0)\Big| \leq
\frac{\lambda L h^{\beta+1}}{(\lfloor \beta\rfloor +1)!}
\sum_{k\in M_{D_{x_0}}} |a_k(x_0)| \leq C_2 \lambda L h^\beta,
\end{eqnarray*}
where the last inequality follows from (\ref{eq:a-coefficients}).
\epr
\par\medskip 
Now we complete the proof of Theorems~\ref{th:upper-bound}.
\par 
First we note that 
because $G\in \sM_2(K)$,
\begin{eqnarray}
 \sum_{i=1}^n H_i&=&\sum_{i=1}^n H(i\delta) \leq \tfrac{1}{\delta} \int_0^{T} H(t) \rd t 
\nonumber
\\
&=&\tfrac{\mu}{\delta}
\int_0^{T}\int_t^\infty [1-G(x)] \rd x\rd t \leq \tfrac{\mu}{\delta}\int_0^\infty x [1-G(x)]\rd x
\leq \tfrac{\mu}{2\delta} K.
\label{eq:sum-H}
\end{eqnarray}
Moreover,  $G\in \sM_2(K)$ implies also that $\frac{1}{\mu}\leq \sqrt{K}$. 
\par 
It can be easily verified that under (\ref{eq:x-0-condition}) and   
(\ref{eq:d-*}) for all $T$ large enough we have 
$T-x_0\geq \kappa T$, and 
$D_{x_0}$ contains at least $\ell +1$ grid points. 
Therefore, by Lemmas~\ref{lem:stoch-error} and~\ref{lem:bias}
and (\ref{eq:sum-H}), the chosen window width $h=h_*$ balances
the upper bounds on the two terms on the right hand side of (\ref{eq:bias-variance}).
The result of Theorem~\ref{th:upper-bound} follows immediately by substitution of $h_*$
in the bounds of   Lemmas~\ref{lem:stoch-error} and~\ref{lem:bias}.
\par\medskip
In order to prove   Theorem~\ref{th:estimation-lambda} 
we note that
the bias--variance decomposition in the problem of estimating $\lambda$
takes the form
\begin{eqnarray}
\Big[\rE_{G,\lambda}  |\hat{\lambda}- \lambda|^2\Big]^{1/2}
\;\leq\; 
\Big\{\rE_{G,\lambda} \Big[\sum_{k\in M_{D_{x_0}}} a_k(x_0) (\hat{R}_k - R_k)\Big]^2\Big\}^{1/2} + 
\Big|\sum_{k\in M_{D_{x_0}}} a_k(x_0) R_k - R^\prime(x_0)\Big|;
\nonumber 
\end{eqnarray}
cf. (\ref{eq:bias-variance}).
The same 
upper bounds on the bias
(Lemma~\ref{lem:bias}) and the variance (Lemma~\ref{lem:stoch-error}) hold.
The upper bound (\ref{eq:upper-bound-arrival-rate}) follows
by 
the special choice of the window width in (\ref{eq:d-*-arrival-rate}). 
\epr 
\subsection{Proof of Theorem~\ref{th:lower-bound}}
The following notation and definitions are used throughout the proof. 
\par 
If $A=\{a_{ij}\}_{i,j=1,\ldots,n}$ 
is an $n\times n$ matrix then 
$\|A\|_2=\sup_{\|x\|_2\leq 1} \|Ax\|_2$ is the spectral norm of $A$,
and $\|A\|_F=(\sum_{i,j=1}^n a_{ij}^2)^{1/2}$ is the Frobenius norm of $A$. 
\par 
Let $v$ be an integrable function on $[-\pi, \pi]$; its  Fourier series is given by 
$v(\omega) = \sum_{j=-\infty}^\infty v_j e^{i\omega j}$, $\omega\in [-\pi, \pi]$, where the 
corresponding Fourier coefficients are 
\[
 v_j=\tfrac{1}{2\pi}\int_{-\pi}^\pi v(\omega) e^{-i\omega j} \rd \omega,\;\;\;j\in\bZ.
\]
For an integrable function $v$ on $[-\pi,\pi]$, let 
$T_n(v)$ stand for the $n\times n$ Toeplitz matrix with the elements
\[
 [T_n(v)]_{j,k}=v_{j-k}=\tfrac{1}{2\pi}\int_{-\pi}^\pi v(\omega) e^{-i(j-k) \omega} \rd \omega,\;\;\;\;\;
j,k=1,\ldots, n.
\]
\subsubsection{Auxiliary results}
The following result  is stated and proved in \citeasnoun{Davies}.
\begin{lemma}\label{lem:first-term}
Let $A$ be an $n\times n$ matrix such that $\|A\|_2<1$; then 
\[
 \big|\log{\rm det}(I+A)-{\rm tr}(A)+\tfrac{1}{2}{\rm tr}(A^2)\big| 
\leq \tfrac{1}{3} \|A\|_2 \|A\|_F^2 (1-\|A\|_2)^{-3}.
\]
\end{lemma}
\par 
In the proof of Theorem~\ref{th:lower-bound}  we use 
properties of Toeplitz's matrices which are presented in the next lemma. 
Some of these statements can be viewed as  ``finite sample'' versions of asymptotic results 
from  \citeasnoun{Davies} and 
\citeasnoun{Dzhaparidze}.
\begin{lemma}\label{lem:Toeplitz}
Let $v, u\in \bL_1([-\pi,\pi])\cap \bL_2([-\pi,\pi])$ be functions with the Fourier coefficients 
$\{v_j\}$ and $\{u_j\}$ respectively. 
\begin{itemize}
 \item[{\rm (i)}] Let $\sup_{\omega\in [-\pi,\pi]} v(\omega)\leq M<\infty$; then 
$\|T_n(v)\|_2 \leq M$.  
\item[{\rm (ii)}] Let $\inf_{\omega\in [\pi,\pi]} v(\omega)\geq m>0$; then $\|T_n^{-1}(v)\|_2 \leq m^{-1}$.
\item[{\rm (iii)}] $\|T_n(v)\|_F^2\leq \frac{n}{2\pi}\int_{-\pi}^\pi |v(\omega)|^2\rd\omega$.
\item[{\rm (iv)}] Suppose that $|v(\omega)|\leq M_1<\infty$, and
$\sum_{j=-\infty}^\infty |j| |u_j|^2\leq M_2<\infty$ for some constants $M_1$ and $M_2$; then 
\[
 \|T_n(vu) - T_n(v) T_n(u)\|_F^2 \leq 4M_1^2 M_2.
\]
\item[{\rm (v)}] Let conditions of (iv) hold, and let 
$\inf_{\omega\in [-\pi,\pi]} v(\omega)\geq m>0$; then
\[
 \|T_n(vu) T_n^{-1}(v)-T_n(u)\|_F^2 \leq 4m^{-2} M_1^2 M_2.
\]
\item[{\rm (vi)}] Let conditions of (iv) and (v) hold; then 
 \begin{eqnarray*}
 {\rm tr}\big\{[T_n^{-1}(v) T_n(vu)]^2\big\} 
\leq \tfrac{n}{\pi} \int_{-\pi}^\pi u^2(\omega) \rd \omega + 8m^{-2}M_1^2 M_2,
\end{eqnarray*}
\end{itemize}
\end{lemma}
\pr 
The statements (i), (ii) and (iii) are standard. See \citeasnoun[p.~64]{Grenander-Szego} for (i) and (ii),
while (iii) is an immediate consequence of Parceval's equality:
\begin{eqnarray*}
 \|T_n(v)\|_F^2 &=& \sum_{j=1-n}^{n-1} (n-|j|) |v_j|^2 \leq n \sum_{j=-\infty}^\infty |v_j|^2=
\tfrac{n}{2\pi}  \int_{-\pi}^\pi |v(\omega)|^2 \rd \omega.
\end{eqnarray*}
\par 
(iv). 
Denote $w(\omega):=v(\omega) u(\omega)$.
By Parceval's equality 
$w_j=\sum_{l=-\infty}^\infty v_l u_{j-l}$, $j\in \bZ$. Therefore 
the $(j,k)$th element of matrix $T_n(w)-T_n(v)T_n(u)$ equals
\begin{eqnarray*}
 \sum_{l=-\infty}^\infty v_l u_{j-k-l} - \sum_{l=1}^n v_{j-l}u_{l-k} = 
\sum_{l=-\infty}^\infty v_l u_{j-k-l} - \sum_{l=j-n}^{j-1} v_{l}u_{j-l-k}
\\
= \sum_{l=-\infty}^{j-n-1} v_l u_{j-k-l} + \sum_{l=j}^\infty v_l u_{j-k-l}.
\end{eqnarray*}
Hence 
\[
\|T_n(vu) - T_n(v) T_n(u)\|_F^2\;\leq\;
 2\sum_{j=1}^n\sum_{k=1}^n \Big|\sum_{l=-\infty}^{j-n-1} v_l u_{j-k-l}\Big|^2\;+\;
2\sum_{j=1}^n\sum_{k=1}^n\Big|\sum_{l=j}^\infty v_l u_{j-k-l}\Big|^2.
\]
Consider the first term; the second term is bounded similarly. 
Let $\Delta$ denote the backward shift operator on 
the space of  two--sided sequences: $(\Delta u)_j=u_{j-1}$, $j\in \bZ$.
For fixed $k$ and $n$ let $u^{(k)}(\cdot)$ be the function on $[-\pi,\pi]$ 
whose Fourier coefficients
are $\{(\Delta^k u)_j {\bf 1}(j\geq n+1), j\in \bZ\}$.
Then with the introduced notation,
\begin{eqnarray*}
\sum_{j=1}^n\sum_{k=1}^n \Big|\sum_{l=-\infty}^{j-n-1} v_l u_{j-k-l}\Big|^2
&\leq & 
\sum_{k=1}^n 
\sum_{j=-\infty}^\infty \Big|\sum_{l=-\infty}^{\infty} v_l 
(\Delta^ku)_{j-l} {\bf 1}\{j-l\geq n+1\} \Big|^2
\\
&=& \tfrac{1}{2\pi}\sum_{k=1}^n \int_{-\pi}^\pi v^2(\omega)[u^{(k)}(\omega)]^2 \rd \omega \;
\leq \; M_1^2 \sum_{k=1}^n \tfrac{1}{2\pi} \int_{-\pi}^\pi [u^{(k)}(\omega)]^2 \rd \omega
\\
&=& M_1^2  \sum_{k=1}^n \sum_{l=-\infty}^\infty |(\Delta^ku)_l {\bf 1}\{l\geq n+1\}|^2 
 \leq M_1^2 \sum_{l=1}^\infty l u_l^2\,, 
\end{eqnarray*}
where the second and third lines follow from  
Parceval's equality and the  premise of the statement.
\par\medskip 
(v). We have 
\begin{eqnarray*}
 \|T_n^{-1}(v) T_n(vu)  - T_n(u)\|_F &=& \big\|T_n^{-1}(v)[T_n(vu)-T_n(v)T_n(u)]\big\|_F 
\\
&\leq& \|T_n^{-1}(v)\|_2
\|T_n(vu)-T_n(v)T_n(u)\|_F.
\end{eqnarray*}
Then the statement follows from (ii) and (iv).
\par\medskip 
(vi). We have
\begin{eqnarray*}
 && {\rm tr}\big\{[T_n^{-1}(v) T_n(vu)]^2 - T_n^2(u)]\big\} 
 \\
 &&\;\;\;\;\leq\; 
\|T_n^{-1}(v)T_n(vu)-T_n(u)\|_F \Big[\|T_n(u)\|_F + \|T_n^{-1}(v) T_n(vu)\|_F\Big].
\end{eqnarray*}
Clearly, 
\[
\|T_n^{-1}(v) T_n(vu)\|_F \;\leq\; \| T_n^{-1}(v)T_n(vu)-T_n(u)\|_F + \|T_n(u)\|_F; 
\]
hence,
by (v) and (iii)
\begin{eqnarray*}
&& {\rm tr}\big\{[T_n^{-1}(v) T_n(vu)]^2 - T_n^2(u)]\big\} 
 \\
&& \;\leq\; 
\|T_n^{-1}(v)T_n(vu)-T_n(u)\|_F^2 + 2\| T_n^{-1}(v)T_n(vu)-T_n(u)\|_F\|T_n(u)\|_F
\\
&& \;\leq\; 
4m^{-2}M_1^2 M_2 + 2m^{-1}M_1 \sqrt{M_2}
\Big[\tfrac{n}{2\pi}\int_{-\pi}^\pi u^2(\omega)\rd \omega\Big]^{1/2}.
\end{eqnarray*}
Therefore using
(iii) we obtain 
\begin{eqnarray*}
 {\rm tr}\big\{[T_n^{-1}(v) T_n(vu)]^2\big\} \leq {\rm tr}\{T_n^2(u)\} + 
4m^{-2}M_1^2 M_2 + 2m^{-1}M_1 \sqrt{M_2}
\Big[\tfrac{n}{2\pi}\int_{-\pi}^\pi u^2(\omega)\rd \omega\Big]^{1/2}
\\
\leq \tfrac{n}{\pi} \int_{-\pi}^\pi u^2(\omega) \rd \omega + 8m^{-2}M_1^2 M_2,
\end{eqnarray*}
as claimed.
\epr 
\subsubsection{Proof of Theorem~\ref{th:lower-bound}} 
The proof is based on standard reduction   
to  a two--point hypotheses testing problem [cf. \citeasnoun[Chapter~2]{tsybakov}].
\par 
Throughout the proof the following notation and conventions are used.
We use symbols $c_0, c_1,\ldots, C_0, C_1,\ldots$ to denote 
positive constants depending on $\beta$, $x_0$, $d$ and $K$ only, 
unless explicitly specified.
\par\medskip
0$^0$. {\em Reduction to a hypotheses testing problem.} Let $x_0$ be fixed, 
$I=[x_0-d, x_0+d]\subset (0,\infty)$,
and 
let $\gamma_0$ and $\gamma_1$ be a pair of covariance functions from $\sC_\beta(L,I, K)$ 
such that 
\begin{equation}\label{eq:a0}
a:=|\theta(\gamma_0)-\theta(\gamma_1)|=|\gamma^\prime_0(x_0)-\gamma^\prime_1(x_0)|>0.
\end{equation}
For arbitrary estimator $\hat{\theta}$ we have 
\begin{eqnarray}
 \cR_{x_0}[\hat{\theta}; \sC_\beta(L,I, K)] \;\geq\; \sup_{\gamma\in \{\gamma_0, \gamma_1\}}
\rE_\gamma |\hat{\theta} - \theta(\gamma)|^2
\;\geq \; \tfrac{1}{4}a^2 \sup_{\gamma\in \{\gamma_0, \gamma_1\}}
\rP_\gamma \big\{|\hat{\theta} - \theta(\gamma)|\geq \tfrac{a}{2}\big\}
\nonumber
\\
\geq \; 
\tfrac{1}{8}a^2 \Big[ \rP_{\gamma_0}\big\{|\hat{\theta} - \theta(\gamma_0)|\geq \tfrac{a}{2}\big\}
+ \rP_{\gamma_1}\big\{|\hat{\theta} - \theta(\gamma_1)|\geq \tfrac{a}{2}\big\}\Big].
\label{eq:reduction}
\end{eqnarray}
\par 
Suppose that on the basis of observations $X^n=(X_1,\ldots, X_n)$ we want to test the hypothesis
$H_0: \theta(\gamma)=\theta(\gamma_0)$ 
against the alternative $H_1: \theta(\gamma)=\theta(\gamma_1)$. Assume that for this purpose
we apply the following {\em minimum distance} testing procedure $\psi(X^n)$: 
given an estimator $\hat{\theta}$ we accept the $i$th hypothesis, $i=0,1$ with 
$\theta(\gamma_i)$ closest to $\hat{\theta}$, i.e., 
$\psi(X^n)={\rm arg}\min_{i=0,1} |\hat{\theta}-\theta(\gamma_i)|$. 
Then, by the triangle inequality, the expression on the right hand side of (\ref{eq:reduction})
is not less than the sum of error probabilities of the minimum distance test:
\begin{eqnarray}\label{eq:low}
 \cR_{x_0}[\hat{\theta}; \sC_{\beta}(L,I,K)] &\geq&
\tfrac{1}{8}a^2 \big[ \rP_{\gamma_0}\{\psi(X^n)=1\}
+ \rP_{\gamma_1}\big\{\psi(X^n)=0\}\big]
\;\geq\;\tfrac{1}{8} a^2 \pi (\rP_{\gamma_0}, \rP_{\gamma_1}), 
\end{eqnarray}
where $\pi (\rP_{\gamma_0}, \rP_{\gamma_1})=\inf_{\varphi} [\rE_{\gamma_0}(1-\varphi)+\rE_{\gamma_1}\varphi]$
is the {\em testing affinity} between $\rP_{\gamma_0}$ and $\rP_{\gamma_1}$; the infimum is taken over all 
tests measurable with respect to the  observation $X^n$.
\par 
Thus the problem is reduced to constructing the worst--case alternatives $\gamma_0$ and $\gamma_1$, and
bounding the testing affinity $\pi(\rP_{\gamma_0}, \rP_{\gamma_1})$.
The last step will be accomplished by bounding from above the 
Kullback--Leibler divergence 
$\cK(\rP_{\gamma_0}, \rP_{\gamma_1})=\rE_{\gamma_1} \log\big [(\rd \rP_{\gamma_1}/\rd \rP_{\gamma_0})(X^n)\big]$ 
between $\rP_{\gamma_0}$ and $\rP_{\gamma_1}$ because 
\begin{equation}\label{eq:low-1}
\pi (\rP_{\gamma_0}, \rP_{\gamma_1}) \geq \tfrac{1}{2}\exp\{-\cK(\rP_{\gamma_0}, \rP_{\gamma_1})\}; 
\end{equation}
see, e.g., \citeasnoun[Theorem~2.4.2]{tsybakov}.
\par\medskip 
1$^0$. {\em Construction of the worst--case alternatives.} 
Let 
$\widehat{\phi}$ be an infinitely differentiable even function with the following properties:
\begin{equation}\label{eq:phi-hat}
 \widehat{\phi}(\omega)=\left\{\begin{array}{ll}
       1, & |\omega|\leq 1,
\\
0, & |\omega|>3/2,                         
                               \end{array}
\right.
\;\;\;\;0\leq \widehat{\phi}(\omega) \leq 1,\;\;
\omega\in [-\tfrac{3}{2}, -1]\cup [1,\tfrac{3}{2}],
\end{equation}
and $\widehat{\phi}$ is monotone on $[-\frac{3}{2},-1]$ and $[1,\frac{3}{2}]$.
Because $\widehat{\phi}$ is an infinitely differentiable function with bounded support, 
the inverse Fourier transform $\phi$ of $\widehat{\phi}$
is a rapidly decreasing infinitely differentiable function: for all $m, k\in \bN$
\begin{equation}\label{eq:fast-decay}
 (it)^m \phi^{(k)}(t) = 
\tfrac{1}{2\pi} \int_{-\infty}^\infty (-i\omega)^k \widehat{\phi}^{(m)}(\omega) e^{-i\omega t}\rd \omega
\;\;
\Rightarrow\;\;
|\phi^{(k)}(t)| \leq C(m,k) |t|^{-m},\;\;\forall t.
\end{equation}
\par 
Let $\ell$ be an even integer number, and let 
\[
\zeta_\ell (t)=(\underbrace{{\bf 1}_{[-1,1]}\ast\cdots\ast 
{\bf 1}_{[-1,1]}}_{\ell})(t),\;\;\;
t\in \bR,
\]
where $\ast$ stands for the convolution on $\bR$.
Put $\zeta(t):= \zeta_\ell(\ell t/(x_0-d))$. Clearly,
${\rm supp}(\zeta)=[-x_0+d, x_0-d]$, and 
the Fourier transform of $\zeta$ is
\begin{equation}\label{eq:zeta-hat}
\widehat{\zeta}(\omega)= 
\int_{-\infty}^\infty \zeta(t) e^{i\omega t}\rd t =
[(x_0-d)/\ell]\bigg[\frac{2\sin\big(\omega(x_0-d)/\ell\big)}{\omega(x_0-d)/\ell}\bigg]^\ell,\;\;\;\omega\in \bR.
\end{equation}
Because $\ell$ is even, function $\widehat{\zeta}$ is non--negative on $\bR$.  
\par
Let 
\[
N=\tfrac{2\pi}{x_0}(N_0+\tfrac{1}{4}), 
\]
where $N_0\geq 1$ is an integer number to be specified, and  
define 
\begin{equation}\label{eq:f-0}
 f_0(\omega) = c_0\delta \widehat{\phi}(\omega\delta/\pi) + c_1 
\big[ \widehat{\zeta}(\omega-N) + \widehat{\zeta}(\omega+N)\big],
\end{equation}
where   $c_0$ and $c_1$ are positive constants. 
We claim that, 
under appropriate choice of 
$c_0$ and $c_1$ function $f_0$
is a spectral density with the corresponding covariance function 
$\gamma_0$ that belongs to $\sC_\beta (c_2L, I, K)$ with preassigned 
$c_2\in (0,1)$.  
\par 
By definition $f_0$ is non--negative and even on $\bR$; hence by Bochner's theorem, 
it is a spectral density. The corresponding covariance function is
\begin{eqnarray*}
 \gamma_0(t) = \tfrac{1}{2\pi}\int_{-\infty}^\infty f_0(\omega)e^{-i\omega t}\rd \omega
= \tfrac{c_0}{2} \phi(\pi t/\delta) +  \tfrac{c_1}{\pi} \zeta(t) \cos(Nt).
\end{eqnarray*}
Because 
${\rm supp}(\zeta)=[-x_0+d, x_0-d]$, 
$\gamma_0(t)=\tfrac{c_0}{2} \phi(\pi t/\delta)$ for $t\in I=[x_0-d, x_0+d]$.
Then 
in view of (\ref{eq:fast-decay}),
\[
 |\gamma_0^{(\beta+1)}(t)| = \tfrac{c_0}{2}  |\phi^{(\beta+1)}(\pi t/\delta)| 
(\pi/\delta)^{(\beta+1)} \leq \tfrac{c_0}{2}  C(\beta+1,\beta+1) |x_0-d|^{-\beta-1}, 
\;\;\;\forall t\in I,
\]
where $C(\cdot, \cdot)$ is a constant appearing in (\ref{eq:fast-decay}).
Choosing $c_0$ and $c_1$ small enough we ensure that  
\mbox{$\gamma_0\in \sC_\beta(c_2L,I,K)$} with  preassigned $0<c_2<1$.
\par\medskip 
Now we proceed with definition of $f_1$ and $\gamma_1$.
With $f_0$ and $\widehat{\phi}$ given by (\ref{eq:f-0}) and (\ref{eq:phi-hat}) respectively
define
\[
 \psi(\omega):=f_0(\omega) \omega \sin(\omega x_0) \big[\widehat{\phi}\big(\tfrac{6x_0}{\pi}(\omega-N)\big) + 
\widehat{\phi}\big(\tfrac{6x_0}{\pi}(\omega+N)\big)\big].
\]
By definition,
\[
 {\rm supp}(\psi) =
\big[-N-\tfrac{\pi}{4x_0}, -N+\tfrac{\pi}{4x_0}\big]\cup \big[N-\tfrac{\pi}{4x_0}, N+\tfrac{\pi}{4x_0}
\big].
\]
For a function $g$ on $\bR$
we put
\[
 B_N(g) := \int_{-\infty}^\infty
g^2(\omega) \sin^2(\omega x_0) \omega^2 \big[\widehat{\phi}\big(\tfrac{6x_0}{\pi}(\omega-N)\big) + 
\widehat{\phi}\big(\tfrac{6x_0}{\pi}(\omega+N)\big)\big]\rd\omega,
\]
and let
\begin{equation}\label{eq:f1}
 f_1(\omega)=f_0(\omega)\big\{1+c_3LN^{-\beta} [B_N(f_0)]^{-1}\psi(\omega)\big\}.
\end{equation}
\par 
Let us verify that $f_1$ is a spectral density.
We have $f_1(\omega)=f_1(-\omega)$ for all $\omega\in \bR$.
To ensure that $f_1$ is non--negative it suffices to require that 
\begin{equation}\label{eq:positivity-f1}
 1\;\geq\; c_3 LN^{-\beta} [B_N(f_0)]^{-1}f_0(\omega)|\omega|,
 \;\;\;\forall \omega \in {\rm supp}(\psi).
\end{equation}
 By definition of $\widehat{\phi}$, $B_N$ and $f_0$,
\begin{eqnarray}
 B_N(f_0) &\geq& 
 2\int_{N-\frac{\pi}{6x_0}}^{N+\frac{\pi}{6x_0}} f_0^2(\omega) \sin^2(\omega x_0) \omega^2\rd \omega 
\nonumber
\\
&\geq& c_4 \int_{N-\frac{\pi}{6x_0}}^{N+\frac{\pi}{6x_0}} [\widehat{\zeta}(\omega-N)]^2 \omega^2\rd\omega 
=  c_4\int_{-\frac{\pi}{6x_0}}^{\frac{\pi}{6x_0}} [\widehat{\zeta}(\omega)]^2 (\omega+N)^2\rd\omega
\geq c_5  N^2,
\label{eq:B-N-1}
\end{eqnarray}
where 
we took into account that  $\sin (\omega x_0)\geq \sqrt{3}/2$ 
whenever $\omega \in [N-\frac{\pi}{6x_0}, N+\frac{\pi}{6x_0}]$.
Moreover, $f_0(\omega)\leq c_6$ for all $\omega\in \bR$; therefore (\ref{eq:positivity-f1}) 
will hold if  $1 \geq c_7LN^{-\beta-1}$. 
This condition will be ensured for large $T$ by our final 
choice of $N$ [cf. (\ref{eq:N-choice})]. 
Thus $f_1$ is a non--negative function on $\bR$, and hence a spectral density.
\par
Let $\gamma_1$ be the covariance function corresponding to $f_1$.
It is evident that choosing $c_3$ small enough we can guarantee that 
$\int_{-\infty}^\infty |\gamma_1(t)|\rd t\leq K$.
It remains to check the smoothness condition. For this purpose we observe that 
\begin{eqnarray*}
&& c_3LN^{-\beta} [B_N(f_0)]^{-1}\int_{-\infty}^\infty |\omega|^{\beta+1} |f_0(\omega)\psi(\omega)|\rd \omega  
\\
&& \;\leq\; 
2c_3LN^{-\beta} [B_N(f_0)]^{-1}\int_{0}^{\infty} 
\omega^{\beta+2} f_0^2 (\omega)|\sin(\omega x_0)|\,
\widehat{\phi} \big(\tfrac{6x_0}{\pi}(\omega-N)\big) \rd\omega 
\\
&&
\;\leq\; 
2c_3 LN^{-\beta} \big(N+\tfrac{\pi}{4x_0}\big)^\beta [B_N(f_0)]^{-1}\int_{0}^{\infty} 
\omega^{2} f_0^2 (\omega)|\sin(\omega x_0)|\,
\widehat{\phi} \big(\tfrac{6x_0}{\pi}(\omega-N)\big) \rd\omega 
\;\leq\; c_8 L,
\end{eqnarray*}
where the last line follows  from 
the fact that  $\widehat{\phi} (\tfrac{6x_0}{\pi}(\cdot - N))$
is supported on $[N-\pi/(4x_0), N+\pi/(4x_0)]$, $\sin (\omega x_0)\geq 1/\sqrt{2}$ on this interval, and 
from the definition of $B_N(f_0)$.
This, together with $\gamma_0\in \sC_\beta(c_2L,I, K)$, means that $\gamma_1\in \sC_{\beta}(L,I, K)$
by choice of constant $c_3$.  
\par\medskip
2$^0$. {\em Distance between the estimated values.} 
We have
\begin{eqnarray}
 a&=&|\gamma_0^\prime(x_0)-\gamma_1^\prime(x_0)|=
\tfrac{1}{\pi} \Big|\int_0^\infty [f_0(\omega)-f_1(\omega)]\omega \sin (\omega x_0) \rd \omega \Big|
\nonumber
\\
&=& \tfrac{c_3}{\pi}LN^{-\beta} [B_N(f_0)]^{-1} \int_0^\infty f_0^2(\omega) \omega^2\sin^2(\omega x_0) 
\widehat{\phi}\big(\tfrac{6x_0}{\pi}(\omega-N)\big)\rd\omega = \tfrac{c_3}{2\pi} L N^{-\beta}, 
\label{eq:a}
\end{eqnarray}
where the last equality follows from definition of $B_N(f_0)$.
\par\medskip 
3$^0$. {\em Spectral densities of the sampled discrete--time process.}  
For generic function $g$ on $\bR$ denote
\[
 \widetilde{g} (\omega) := \tfrac{1}{\delta} \sum_{j=-\infty}^\infty g\big(\tfrac{\omega + 2\pi j}{\delta}\big),\;\;\;
\omega \in (-\pi, \pi].
\]
Under $\rP_{\gamma_0}$ and $\rP_{\gamma_1}$, 
the spectral densities of the discrete--time process $\{X(k\delta), k\in\bZ\}$,
are $\widetilde{f}_0$ and $\widetilde{f}_1$  respectively;
see, e.g., \citeasnoun{Grenander-Rosenblatt}.
By~(\ref{eq:f1})
\begin{eqnarray}\label{eq:f-1-tilde}
 \widetilde{f}_1(\omega) &=& \widetilde{f}_0(\omega) + c_3 LN^{-\beta} [B_N(f_0)]^{-1}\,
\tfrac{1}{\delta}
\sum_{j=-\infty}^\infty f_0\big(\tfrac{\omega+2\pi j}{\delta}\big) \psi\big(\tfrac{\omega+2\pi j}{\delta}\big).
\end{eqnarray}
In what follows we require that 
\begin{equation}\label{eq:N-delta}
(N+\tfrac{\pi}{4x_0})\delta \leq \pi;
\end{equation}
this condition will be verified by our choice of~$N$.  
Under this condition, since function $\psi$
is supported on 
$[N-\pi/(4x_0), N+\pi/(4x_0)]\cup [-N-\pi/(4x_0), -N+\pi/(4x_0)]$, 
the sum in (\ref{eq:f-1-tilde}) 
contains only one non--vanishing term corresponding to $j=0$. Thus, 
\begin{eqnarray*}
 \widetilde{f}_1(\omega) = \widetilde{f}_0(\omega) + c_3\delta^{-1} LN^{-\beta} [B_N(f_0)]^{-1}
f_0(\omega/\delta) \psi(\omega/\delta) 
= \widetilde{f}_0(\omega)[1 + g(\omega)],
\end{eqnarray*}
where we have denoted
\begin{equation}\label{eq:g}
 g(\omega) := c_3\delta^{-1}LN^{-\beta} 
\frac{ f_0^2(\omega/\delta)(\omega/\delta)}{B_N(f_0)\widetilde{f}_0(\omega)}
\sin(\omega x_0/\delta)
\Big[\widehat{\phi}\big(\tfrac{6x_0}{\pi}(\tfrac{\omega}{\delta}-N)\big)+ 
 \widehat{\phi}\big(\tfrac{6x_0}{\pi}(\tfrac{\omega}{\delta}+N)\big)\Big].
\end{equation}
\par 
The next lemma summarizes some useful properties of function $g$.
\begin{lemma}\label{lem:properties-g}
The following statement holds:
\begin{equation}\label{eq:g-l2}
 \int_{-\pi}^\pi |g(\omega)|^2\rd \omega \;\leq\; c_3^2L^2 N^{-2\beta} [B_N(f_0)]^{-1}\delta. 
\end{equation}
In addition, if $\{g_j\}$ are the Fourier coefficients of the function $g$ then
\begin{equation}\label{eq:j-g}
 \sum_{j=-\infty}^\infty |j| |g_j|^2 \;\leq\;  c_{9} L^2N^{-2\beta} [B_N(f_0)]^{-3/2} 
\Big\{ [B_N(f_0)]^{1/2}+[B_N(f^\prime_0)]^{1/2} + [B_N(f^\sharp_0)]^{1/2}\Big\},
\end{equation}
where $f_0^\sharp(\omega):=\delta \widetilde{f^\prime_0}(\omega)
=\sum_{j=-\infty}^\infty f_0^\prime (\omega +\frac{2\pi j}{\delta})$.
\end{lemma}
\pr Since $\widetilde{f}_0(\omega)\geq \delta^{-1} f_0(\omega/\delta)$, $\forall \omega\in [-\pi, \pi]$
 we obtain 
\begin{eqnarray*}
 && \int_{-\pi}^\pi |g(\omega)|^2 \rd \omega
\\
 && \;\;\leq\;
c_3^2L^2N^{-2\beta}[B_N(f_0)]^{-2}\int_{-\pi}^\pi f_0^2(\omega/\delta)
(\omega/\delta)^2 \sin^2(\omega x_0/\delta) 
\Big[\widehat{\phi}\big(\tfrac{6x_0}{\pi}(\tfrac{\omega}{\delta}-N)\big)+ 
 \widehat{\phi}\big(\tfrac{6x_0}{\pi}(\tfrac{\omega}{\delta}+N)\big)\Big]^2\rd\omega
\\
&&\;\;=\;
c_3^2L^2N^{-2\beta}[B_N(f_0)]^{-2}\delta \int_{-\pi/\delta}^{\pi/\delta} f_0^2(\omega)
\omega^2 \sin^2(\omega x_0) 
\Big[\widehat{\phi}\big(\tfrac{6x_0}{\pi}(\omega-N)\big)+ 
 \widehat{\phi}\big(\tfrac{6x_0}{\pi}(\omega+N)\big)\Big]^2\rd\omega
\\
&&\;\;\leq\; c_3^2 L^2N^{-2\beta}[B_N(f_0)]^{-1}\delta, 
\end{eqnarray*}
as claimed. The last inequality follows by definitions of $\widehat{\phi}$ and $B_N(\cdot)$. 
\par 
Now we prove the second statement of the lemma.  
Write for brevity $A=c_3LN^{-\beta}[B_N(f_0)]^{-1}$ and note that 
$g(\omega)=g_0(\omega/\delta)$,
 where 
\begin{equation}\label{eq:tilde-g}
 g_0(\omega) := A 
\Big[\sum_{j=-\infty}^\infty f_0\big(\omega + \tfrac{2\pi j}{\delta}\big)\Big]^{-1}
f_0^2(\omega) \omega \sin (\omega x_0)\big[\widehat{\phi}\big(\tfrac{6x_0}{\pi}(\omega-N)\big)+
\widehat{\phi}\big(\tfrac{6x_0}{\pi}(\omega+N)\big)\big];
\end{equation}
see (\ref{eq:g}).
By the Cauchy--Schwarz inequality and Parceval's equality
\begin{eqnarray*}
&&\sum_{j=-\infty}^\infty |j| |g_j|^2\leq 
\Big(\sum_{j=-\infty}^\infty |g_j|^2\Big)^{1/2}\Big(\sum_{j=-\infty}^\infty j^2 |g_j|^2\Big)^{1/2}
\\
&& =
\tfrac{1}{2\pi}\Big(\int_{-\pi}^\pi |g(\omega)|^2\rd \omega \Big)^{1/2}
\Big(\int_{-\pi}^\pi |g^\prime(\omega)|^2\rd \omega \Big)^{1/2}\;
\leq c_3LN^{-\beta}[B_N(f_0)]^{-1/2}  
\Big(\int_{-\pi/\delta}^{\pi/\delta} |g_0^\prime(\omega)|^2\rd \omega\Big)^{1/2},
\end{eqnarray*}
where in the last step we used (\ref{eq:g-l2}). 
We proceed with bounding the integral on the right hand side. 
\par
It follows from (\ref{eq:tilde-g}) that
 $g_0^\prime (\omega) = \sum_{m=1}^5 J_m(\omega)$, 
where 
\begin{eqnarray*}
&& J_1(\omega) = A 
\Big[\sum_{j=-\infty}^\infty f_0\big(\omega + \tfrac{2\pi j}{\delta}\big)\Big]^{-1}
2f_0(\omega) f_0^\prime(\omega) \omega \sin (\omega x_0)
\big[\widehat{\phi}\big(\tfrac{6x_0}{\pi}(\omega-N)\big)+
\widehat{\phi}\big(\tfrac{6x_0}{\pi}(\omega+N)\big)\big],
\\
&&J_2(\omega) = A 
\Big[\sum_{j=-\infty}^\infty f_0\big(\omega + \tfrac{2\pi j}{\delta}\big)\Big]^{-1}
f_0^2(\omega) \sin (\omega x_0)
\big[\widehat{\phi}\big(\tfrac{6x_0}{\pi}(\omega-N)\big)+
\widehat{\phi}\big(\tfrac{6x_0}{\pi}(\omega+N)\big)\big],
\\
&& J_3(\omega) = A 
\Big[\sum_{j=-\infty}^\infty f_0\big(\omega + \tfrac{2\pi j}{\delta}\big)\Big]^{-1}
f_0^2(\omega) \omega x_0\cos (\omega x_0)
\big[\widehat{\phi}\big(\tfrac{6x_0}{\pi}(\omega-N)\big)+
\widehat{\phi}\big(\tfrac{6x_0}{\pi}(\omega+N)\big)\big],
\\
&& J_4(\omega) = A 
\Big[\sum_{j=-\infty}^\infty f_0\big(\omega + \tfrac{2\pi j}{\delta}\big)\Big]^{-1}
f_0^2(\omega) \omega \sin (\omega x_0)\tfrac{6x_0}{\pi}
\big[\widehat{\phi}^\prime\big(\tfrac{6x_0}{\pi}(\omega-N)\big)+
\widehat{\phi}^\prime\big(\tfrac{6x_0}{\pi}(\omega+N)\big)\big],
\\
&& J_5(\omega) = - \frac{A\, f_0^2(\omega) \omega \sin (\omega x_0)}{ 
[\sum_{j=-\infty}^\infty f_0\big(\omega + 2\pi j/\delta\big)]^{2}}\;
\big[\widehat{\phi}\big(\tfrac{6x_0}{\pi}(\omega-N)\big)+
\widehat{\phi}\big(\tfrac{6x_0}{\pi}(\omega+N)\big)\big]
\sum_{j=-\infty}^\infty f_0^\prime\big(\omega+
\tfrac{2\pi j}{\delta}\big).
\end{eqnarray*}
Since $\sum_{j=-\infty}^\infty f_0(\omega+2\pi j/\delta) \geq f_0(\omega)$,
\begin{eqnarray*}
 \int_{-\pi/\delta}^{\pi/\delta}|J_1(\omega)|^2\rd\omega \leq 
4A^2 B_N(f_0^\prime) = 4 c^2_3L^2 N^{-2\beta}[B_N(f_0)]^{-2} B_N(f_0^\prime).
\end{eqnarray*}
Similarly we obtain the following bounds
\begin{eqnarray*}
\int_{-\pi/\delta}^{\pi/\delta} |J_2(\omega)|^2\rd\omega 
&\leq&  A^2 \int_{-\pi/\delta}^{\pi/\delta}
f_0^2(\omega) \sin^2(\omega x_0)\big[\widehat{\phi}\big(\tfrac{6x_0}{\pi}(\omega-N)\big)+
\widehat{\phi}\big(\tfrac{6x_0}{\pi}(\omega+N)\big)\big]^2\rd\omega 
\\
&\leq& A^2 (N-\tfrac{\pi}{4x_0})^{-2} \int_{-\pi/\delta}^{\pi/\delta}
f_0^2(\omega) \omega^2 \sin^2(\omega x_0)\big[\widehat{\phi}\big(\tfrac{6x_0}{\pi}(\omega-N)\big)+
\widehat{\phi}\big(\tfrac{6x_0}{\pi}(\omega+N)\big)\big]\rd\omega 
\\
&\leq& c L^2 N^{-2\beta-2} [B_N(f_0)]^{-1},
\\*[2mm]
\int_{-\pi/\delta}^{\pi/\delta} |J_k(\omega)|^2\rd \omega &\leq& c L^2N^{-2\beta} [B_N(f_0)]^{-1},\;\;\;k=3,4,
\\
\int_{-\pi/\delta}^{\pi/\delta} |J_5(\omega)|^2\rd \omega &\leq& c L^2 N^{-2\beta}[B_N(f_0)]^{-2} 
B_N(f^\sharp_0),
\end{eqnarray*}
where $f^\sharp_0(\omega)=\sum_{j=-\infty}^\infty f^\prime_0 (\omega + \tfrac{2\pi j}{\delta})$.
Combining these results we come to the statement (\ref{eq:j-g}).
\epr
\par\medskip 
4$^0$. {\em The Kullback--Leibler divergence.} 
Now we proceed with bounding the Kullback--Leibler 
divergence between the probability measures $\rP_{\gamma_0}$
and $\rP_{\gamma_1}$ generated by observation $X^n$ under hypotheses $H_0$ and $H_1$.
Under $H_0$  the distribution of observation $X^n$
is multivariate normal with zero mean and 
covariance matrix $\Sigma_0=T_n(\widetilde{f}_0)$, while under $H_1$
 the distribution  is the multivariate normal $\rP_{\gamma_1}$ with zero mean 
and covariance matrix
$\Sigma_1=T_n(\widetilde{f}_1)= T_n(\widetilde{f}_0) + T_n(\widetilde{f}_0g)$.
The Kullback--Leibler divergence between these multivariate normal distributions is 
\begin{eqnarray*}
 \cK(\rP_{\gamma_0}, \rP_{\gamma_1}) = 
\rE_{\gamma_1}\log \Big[\tfrac{\rd \rP_{\gamma_1}}{\rd \rP_{\gamma_0}}(X^n)\Big] = \tfrac{1}{2}
\log \frac{{\rm det}(\Sigma_0)}{{\rm det}(\Sigma_1)} - \tfrac{1}{2}n + \tfrac{1}{2}{\rm tr}
(\Sigma_0^{-1}\Sigma_1).
\end{eqnarray*}
Put for brevity $V=T(\widetilde{f}_1)-T_n(\widetilde{f}_0)=T_n(\widetilde{f}_0g)$; then 
\begin{eqnarray}
\cK(\rP_{\gamma_0}, \rP_{\gamma_1}) = - \tfrac{1}{2}\log {\rm det}(\Sigma_0^{-1}\Sigma_1) + \tfrac{1}{2}
{\rm tr}(\Sigma_0^{-1}\Sigma_1 - I) 
\;=\; -\tfrac{1}{2} \log {\rm det}(I+\Sigma_0^{-1}V) + \tfrac{1}{2}{\rm tr}(\Sigma_0^{-1}V)
\nonumber
\\
= -\tfrac{1}{2} w(\Sigma_0^{-1}V) + \tfrac{1}{4} {\rm tr}\{(\Sigma_0^{-1}V)^2\},
\label{eq:K-L-1}
\end{eqnarray}
where $w(A) := \log {\rm det} (I+A) - {\rm tr}(A) + \tfrac{1}{2} {\rm tr}\{A^2\}$.
Our current goal is to bound the two terms on the right hand side on (\ref{eq:K-L-1}).
\par \medskip 
First we note
the following upper and lower bounds on $\widetilde{f}_0$.
It follows from the definition of $f_0$ and $\widehat{\phi}$ that
\[
 \widetilde{f}_0(\omega)=\tfrac{1}{\delta}\sum_{j=-\infty}^\infty
f_0\big(\tfrac{\omega+2\pi j}{\delta}\big) \geq c_0  \sum_{j=-\infty}^\infty \widehat{\phi}
\big(\tfrac{\omega+2\pi j}{\pi}\big)\geq c_0 
\inf_{\omega\in [-\pi,\pi]} \widehat{\phi}(\omega/\pi) =c_0.
\]
On the other hand, 
\begin{eqnarray*}
 \widetilde{f}_0(\omega) &\leq& c_{10} + c_1  \delta^{-1}\sum_{j=-\infty}^\infty 
\big[\widehat{\zeta}\big(\tfrac{\omega+2\pi j}{\delta}-N\big) + 
\widehat{\zeta}\big(\tfrac{\omega+2\pi j}{\delta}+N\big)\big]
\\
&\leq & c_{10} + c_1 \delta^{-1} 2^\ell (x_0-d)/\ell  +  
c_1\delta^{-1}\sum_{\substack{j=-\infty\\ j\ne 0}}^\infty 
\big[\widehat{\zeta}\big(\tfrac{\omega+2\pi j}{\delta}-N\big) + 
\widehat{\zeta}\big(\tfrac{\omega+2\pi j}{\delta}+N\big)\big]
\\
&\leq& c_{10} + c_{11}\delta^{-1} + c_{12} \delta^{\ell-1}\;\leq\;
c_{13} \delta^{-1};
\end{eqnarray*}
here the third inequality follows 
from (\ref{eq:zeta-hat}) and (\ref{eq:N-delta}).  
Thus we have shown that 
\begin{equation}\label{eq:f0-bounds}
0<c_0 \leq \inf_{\omega\in [-\pi, \pi]} \widetilde{f}_0(\omega) \leq  
\sup_{\omega\in [-\pi, \pi]} \widetilde{f}_0(\omega)\leq c_{13}\delta^{-1}.
\end{equation}
\par 
We are in a position to bound $\cK(\rP_{\gamma_0}, \rP_{\gamma_1})$ in (\ref{eq:K-L-1})
from above.
The statement (vi) of Lemma~\ref{lem:Toeplitz} together with (\ref{eq:g-l2}) and (\ref{eq:j-g})
implies
\begin{eqnarray}
&& {\rm tr}\{(\Sigma_0^{-1}V)^2\} 
= {\rm tr} \big\{ [T_n^{-1}(\widetilde{f}_0) T_n(\widetilde{f}_0g)]^2\big\}
\nonumber
\\
&& \leq  c_{14} L^{2}N^{-2\beta} [B_N(f_0)]^{-1}
\Big\{n\delta
+  \delta^{-2} [B_N(f_0)]^{-1/2} 
\Big( [B_N(f^\prime_0)]^{1/2} + [B_N(f^\sharp_0)]^{1/2}\Big)\Big\}.
\label{eq:first-term}
\end{eqnarray}
This yields the upper bound on the second term on the right hand side of (\ref{eq:K-L-1}).
\par 
Now we proceed with bounding 
$w(\Sigma_0^{-1}V)$. Note that  
\begin{equation}\label{eq:sigma-0-v}
\|\Sigma_0^{-1}V\|_2=
\|T_n^{-1}(\widetilde{f}_0) T_n(\widetilde{f}_0g)\|_2 \leq \|T_n^{-1}(\widetilde{f}_0)\|_2 
\|T_n(\widetilde{f}_0g)\|_2.
\end{equation}
In view of the lower bound in (\ref{eq:f0-bounds}) and by Lemma~\ref{lem:Toeplitz}(ii),
$\|T_n^{-1}(\widetilde{f}_0)\|_2\leq c_0^{-1}$.
Using the definition of  $g$ [see (\ref{eq:g})] we obtain 
\begin{eqnarray*}
 \widetilde{f}_0(\omega) g(\omega) &=& c_3\delta^{-1}LN^{-\beta} [B_N(f_0)]^{-1}\,
f_0^2(\omega/\delta)(\omega/\delta)
\sin(\omega x_0/\delta)
\Big[\widehat{\phi}\big(\tfrac{6x_0}{\pi}(\tfrac{\omega}{\delta}-N)\big)+ 
 \widehat{\phi}\big(\tfrac{6x_0}{\pi}(\tfrac{\omega}{\delta}+N)\big)\Big]
\\
&\leq& c_3\delta^{-1} LN^{-\beta} 
[B_N(f_0)]^{-1} f_0^2(\omega/\delta) 
(\omega/\delta) \;\leq\; c_{15}\delta^{-1}LN^{-\beta+1} [B_N(f_0)]^{-1},
\end{eqnarray*}
where in the last inequality we took into account that
\[
{\rm supp}(g)=\big[\delta( N-\tfrac{\pi}{4x_0}), \delta (N+ \tfrac{\pi}{4x_0})\big]
\cup \big[-\delta(N+\tfrac{\pi}{4x_0}), -\delta (N-\tfrac{\pi}{4x_0})\big],
\]
and the definition of $f_0$  [see (\ref{eq:f-0})]. 
Then it follows from (\ref{eq:sigma-0-v}) and Lemma~\ref{lem:Toeplitz}(i) that 
\[
\|T_n^{-1}(\widetilde{f}_0) T_n(\widetilde{f}_0g)\|_2\leq 
c_{15} \delta^{-1} LN^{-\beta+1} [B_N(f_0)]^{-1}.
\]
Let us require that the choice of $N$ be such that
\begin{equation}\label{eq:N-condition}
c_{15} \delta^{-1} LN^{-\beta+1} [B_N(f_0)]^{-1}\leq 1/2.
\end{equation}
Then  Lemma~\ref{lem:first-term} is applicable, and 
\[
|w(\Sigma_0^{-1}V)| = |w(T_n^{-1}(\widetilde{f}_0) T_n(\widetilde{f}_0g))| \leq 
\tfrac{4}{3}\|T_n^{-1}(\widetilde{f}_0) T_n(\widetilde{f}_0g)\|_F^2.
\]
Since 
$\|T_n^{-1}(\widetilde{f}_0) T_n(\widetilde{f}_0g)\|_F^2 = {\rm tr}\big\{
[T_n^{-1}(\widetilde{f}_0) T_n(\widetilde{f}_0g)]^2\big\}$, we obtain from 
(\ref{eq:N-condition}), (\ref{eq:first-term}) and (\ref{eq:K-L-1}) that 
\begin{equation}\label{eq:K-L-final}
 \cK(\rP_{\gamma_1}, \rP_{\gamma_0}) \leq 
c_{16} L^{2}N^{-2\beta} [B_N(f_0)]^{-1}
\Big\{n\delta
+  \delta^{-2} [B_N(f_0)]^{-1/2} 
\big( [B_N(f^\prime_0)]^{1/2} + [B_N(f^\sharp_0)]^{1/2}\big)\Big\}.
\end{equation}
\par 
Recall that $B_N(f_0)\geq c_5 N^2$; see (\ref{eq:B-N-1}). Moreover, 
\begin{eqnarray*}
 B_N(f_0^\prime) &\leq & 2\int_{N-\frac{\pi}{4x_0}}^{N+\frac{\pi}{4x_0}} 
[f_0^\prime(\omega)]^2\omega^2\rd\omega 
\leq c_{17} 
\int_{-\frac{\pi}{4x_0}}^{\frac{\pi}{4x_0}} 
[\widehat{\zeta}^\prime(\omega)]^2 (\omega+N)^2\rd\omega 
\leq c_{18} N^2,
\end{eqnarray*}
where the last inequality is obtained from the definition of $\zeta(\cdot)$.
By similar argument one can show that   $B_N(f_0^\sharp)\leq c_{19}N^2$, provided that $\ell$ 
is large enough. 
Combining these inequalities with (\ref{eq:K-L-final}) we finally obtain
\[
 \cK(\rP_{\gamma_1}, \rP_{\gamma_0}) \leq c_{20} L^2 N^{-2\beta-2} \big(n\delta + \delta^{-2}\big)=
c_{20}L^2 N^{-2\beta-2}(T+\delta^{-2}).
\]
\par\medskip
5$^0$. {\em Choice of $N$ and proof completion.}
Pick integer $N_0$ such that for some constant $c_{22}$
\begin{equation}
\label{eq:N-choice}
 N=N_*=c_{21} (L^2T)^{1/(2\beta+2)}~.
\end{equation}
With this choice
under condition (\ref{eq:T-range})
we have 
\[
 (N_*+\tfrac{\pi}{4x_0})\delta \leq c_{21}(L^2T)^{1/(2\beta+2)} \delta + \pi (4x_0)^{-1}\delta \leq 
c_{21} C_2 + \pi (4x_0)^{-1}\delta  
\]
which is less than $\pi$ by choice of $c_{21}$ and as $\delta\to 0$.
Thus (\ref{eq:N-delta}) holds for all $\delta$ small. 
\par 
Moreover,  in view of (\ref{eq:B-N-1})
\[
 c_{15} \delta^{-1} L N_*^{-\beta+1} [B_{N_*}(f_0)]^{-1} \leq c_{15} c_5 \delta^{-1} LN_*^{-\beta-1} =
c_{15}c_5 c_{21}^{-\beta-1} \delta^{-1} T^{-1/2}\leq c_{15}c_5 c_{21}^{-\beta-1} C_1^{-1},
\]
where in the last inequality we have used (\ref{eq:T-range}). 
The left hand side of the last inequality is less than $1/2$ for $C_1$ large enough.
Thus (\ref{eq:N-condition}) is fulfilled. 
\par 
Finally, if $N=N_*$ then in view of (\ref{eq:T-range}),
$\cK(\rP_{\gamma_0}, \rP_{\gamma_1}) \leq c_{22}$.  Therefore the theorem statement 
follows from (\ref{eq:a}), (\ref{eq:low}) and (\ref{eq:low-1}).
The proof is completed. 
\epr

\bibliographystyle{agsm}

\end{document}